\documentclass[11pt]{article}
\usepackage[margin=1in]{geometry}
\usepackage{amsmath,amssymb,fancyhdr,amsfonts,amsthm,verbatim,graphicx}
\usepackage{latexsym}
\usepackage{exscale}
\usepackage{esvect}
\usepackage{pgf,tikz}
\usepackage{mathrsfs}
\usetikzlibrary{arrows}
\usepackage{cancel}
\usepackage{float}
\usepackage{bm} 
\usepackage{caption}
\usepackage{bigints}
\usepackage{soul}
\usepackage[shortlabels]{enumitem}
\usepackage{ulem}
\usetikzlibrary{trees}

\theoremstyle{definition}
\newtheorem{theorem}{Theorem}[section]

\newtheorem{lemma}[theorem]{Lemma}

\newtheorem{definition}[theorem]{Definition}
\newtheorem{example}[theorem]{Example}

\numberwithin{equation}{section}

\def \beq{\begin{equation}}
	\def \eeq{\end{equation}}
\def \lab{\label}
\renewcommand{\rq}[1]{(\ref{#1})}

\renewcommand{\rq}[1]{(\ref{#1})}

\newtheorem{prop}{Proposition}
\newtheorem{thm}{Theorem}
\newtheorem{cor}{Corollary}
\newcommand{\bD}{\Bbb D}
\newcommand{\bR}{{ \mathbb R  }}
\newcommand{\bC}{\Bbb C}
\newcommand{\bZ}{\Bbb Z}
\newcommand{\bN}{\Bbb N}
\newcommand{\bK}{\Bbb K}
\newcommand{\ov}{\bar v}
\newcommand{\E}{\mathcal{E}}
\newcommand{\F}{\mathcal{F}}
\newcommand{\R}{\mathcal{R}}
\newcommand{\de}{\dot e}
\newcommand{\hg}{\hat g}
\newcommand{\hu}{\hat u}
\newcommand{\V}{\mathop{\rm V}}
\newcommand{\ev}{\mbox{eigenvalues}}
\newcommand{\tphi}{{\tilde \phi}}
\newcommand{\cC}{{\cal C}}
\newcommand{\Res}{\mbox{$Res({\lambda})$}}
\newcommand{\la}{\mbox{$\lambda$}}
\newcommand{\smooth}{\mbox{$\Psi^{-\infty}$}}
\newcommand{\ra}{\mbox{$\mapsto $}}
\newcommand{\Co}{C_0^{\infty}}
\newcommand{\pa }{\partial }
\newcommand{\Da}{\Delta_{a,r}}
\newcommand{\og}{{\overline g}}
\newcommand{\Pf}{\noindent {\it Proof:}\ }
\newcommand{\ep}{\epsilon}
\newcommand{\frf}{(\frac{f'}{f})}
\newcommand{\tb}{\tilde{B}}
\newcommand{\ty}{\tilde{y}}
\newcommand{\N}{{\cal N}}
\newcommand{\al}{\alpha }
\newcommand{\Om}{\Omega}

\newcommand{\f}{\varphi }

\newcommand{\Ga}{\Gamma }

\newcommand{\mH}{{\mathcal H}}

\newcommand{\om}{{\omega}}

\newcommand{\ka}{{\kappa}}

\def\<{\langle} \def\>{\rangle}

\newcommand{\tf}{\tilde{f}}

\newcommand{\hH}{\hat{H}}

\newcommand{\je}[1]{{\color{black}{#1}}}

\title{Exact controllability for wave equation on general quantum graphs with non-smooth controls}
\author{Sergei A. Avdonin and Julian K. Edward}
\date{January 1, 2023}

\begin{document}
	\maketitle
	
\noindent	{\bf Abstract}

{In this paper we study the exact controllability problem for the wave equation on a finite metric graph with the Kirchhoff-Neumann matching conditions.  Among all vertices and edges we choose certain active vertices and edges, and give a constructive proof that the wave equation on the graph is exactly controllable if $H^1(0,T)'$ Neumann controllers are placed at the active vertices and $L^2(0,T)$ Dirichlet controllers are placed at the active edges. The proofs for the shape and velocity controllability are purely dynamical,  while the proof for the exact controllability utilizes both dynamical and moment method approaches.  The control time for this construction is determined by the chosen orientation and path decomposition of the graph.  }

	\section{Introduction}
	Controllability properties of the wave equation is a fundamental topic in the control theory for partial differential equations. Many powerful methods were used to prove controllability of the wave equation in various spatial domains under the action of various types of controls (see, e.g. \cite{LT, Li, Rus, Z07} and references therein).  Control problems for the wave equation on graphs have important applications in science and engineering and were studied in many papers (see  the monographs \cite{AI, DZ,  LLS12};  the surveys \cite{A2, Z13};
	and references therein). They also have deep connection with inverse problem on graphs, see, e.g. \cite{AK, B} for tree graph examples.  In this paper we consider the
	exact controllability problem for the equation	$ u_{tt}-u_{xx}+q(x)u=0 $
	on metric graphs with the Kirchhoff-Neumann matching conditions.
	For graphs without cycles (trees) this problem was studied (in slightly different but essentially equivalent  forms) in, e.g.  \cite{AM08, B, DZ,LLS, LLS12}. 
	It was proved that the system is exactly controllable
	if  Dirichlet controllers act at all or at all but
	one of the boundary vertices.  In \cite{AZ} a constructive proof for controllability of the wave equation on a  tree graph was proposed, in which the  graph was represented as a union of disjoint paths with each path starting from a controlled boundary vertex. The sharp controllability time was also indicated in that paper.
	
	Until recently, little was known about controllability of the wave equation (or any other partial differential equation) on graphs with cycles. It was only proved that the wave equation on graphs with cycles
	is never exactly boundary controllable \cite{AI}; it may be spectrally (boundary) controllable, but this property is very unstable with respect to small perturbations of the system parameters (see, e.g. \cite{ABP,ABI}). 
	Historically, the first work on controllability of the wave equation on graphs was the paper by S. Rolewicz \cite{Rol}. He considered the system  with Dirichlet controls at all vertices, i.e. the system of strings connected by common controls at the internal vertices.
		Complete solution of the problem of exact and spectral controllability for such systems was obtained in \cite[Sec. VII.1]{AI}.
	
	An important step forward recently was the work by \cite{AZ2}, who gave
	a constructive proof for the controllability of the wave equation on a general graph with the Kirchhoff-Neumann matching conditions.  Their construction generalized the idea of controllability on trees in  \cite{AZ}.  However,
	to prove controllability of systems on graphs with cycles requires not only boundary but also internal controls, as was proposed in \cite{A5}.
	Given an undirected graph, they first gave it an acyclic orientation, and then constructed an active set of vertices and edge directions based on the chosen orientation, such that the active vertices include all source vertices specific to the orientation; the active edge directions include all but one of the outgoing directions at every vertex. Then Neumann controllers are placed at the active vertices, and Dirichlet controllers are placed at the active edge directions.  The shape and velocity controllability are proved using a constructive dynamical method.  The exact  controllability was then proved by combining the shape and velocity control functions via the  method of moments.
	
	One important restriction made in \cite{AZ2} is that the
	controls had relatively high regularity: 
	Dirichlet controls were assumed to be in a subset of $H^1(0,T)$, while the Neumann controls were in $L^2(0,T),$ where $T$ is the appropriate controllability time.  It is natural to consider a rougher set of controls. In this paper, we prove an exact controllability assuming $L^2$-Dirichlet controls and $H^1(0,T)'$-Neumann controls, where $H^1(0,T)'$ is the dual of $H^1(0,T).$
	To the best of our knowledge, for such controls it is not even known the appropriate state spaces for which the initial boundary value problem (IBVP) is well posed. In addition, there are significant technical obstructions to implementing the local, dynamical arguments used in \cite{AZ2} for proving velocity controllability. In particular, we are unaware of any natural way to define the restriction $H^1(0,T)'\mapsto H^1(a,b)'$ where $(a,b)$ is a proper subset of $(0,T)$. To avoid these issue, we instead prove velocity control in a  {\it more}
	regular category, and then use the identity 
	$$\frac{\pa^2}{\pa t^2}u^{\bf f}(x,t)=u^{{\bf f}''}(x,t)
	$$
	to obtain velocity control in the rougher category. Here ${\bf f}$ is the vector of controls, and $u^{{\bf f}}$ is the associated solution to the IBVP. 
	This argument can be compared with the work of \cite{EZ}, where controllability in {\it less} regular spaces is used to prove controllability in {\it more} regular spaces. 
	
	The paper is organized as follows.  Section \ref{pre1} introduces the quantum graph and formulates the control problem,  then discusses well-posedness. 
	In Section \ref{sv},
	the orientation of the graph and the placement of the controls is presented, and then shape and velocity control results are proven. For readability we assume in this section that $q=0$; the necessary modifications for $q\neq 0$ are given in the Appendix.  In    Section \ref{Econt}, we formulate the shape and velocity controllability results in terms of moment problems, 
	and then prove the solvability of the moment problem associated to exact controllability. 
	
	
	\section{Well-posedness}
	\subsection{Preliminaries}\label{pre1}
	{
		Let $\Omega(V,E)$ be a finite,  undirected, connected graph, where $V,E$ are the sets of vertices and edges of $\Omega$ respectively.  Each edge in $E$ is associated with two vertices in $V$ called its endpoints.  Here we require both $V = \{v_i: i \in I\}$ and $E= \{e_j: j \in J\}$ to be nonempty; and the two endpoints of an edge to be distinct. That is, our graphs contain no loops;  we will use this assumption in Section \ref{sv} when introducing acyclic orientations on graphs. This assumption is not restrictive.  Since we admit multiple edges between two vertices, a loop can be considered as two edges between two vertices. In this paper it is convenient to consider $I$ and $J$ as disjoint subsets of $\bN.$ 
		
		An undirected edge $e_j$ between $v_i$ and $v_k$ is denoted as $e_j(v_i,v_k)$.  The set of indices of the edges  incident to $v_i$ is denoted by $J(v_i)$.  The set, $ \Gamma=\{v_i \in  V: |J(v_i)|=1\}, $  plays the role of the graph boundary.  Here we use $|\cdot|$ to denote the cardinality of a set of features on a graph.  
		
		We recall that a graph is called a {\bf metric graph} if every edge $ e_j \in E$  is identified with an interval $(0,l_j) $ of the real line with a positive length $l_j$. 
		The graph $ \Omega $ determines naturally the Hilbert space of square
		integrable functions $ \mathcal{H}$ on $\Omega $ such that
		$y_j:=y|_{e_j} \in L^2(0,l_j)$ for every $j \in J.$
		{We define the space  $\mathcal{H}^1$
		of  continuous functions $y$ on $\Omega$ such that
		$y_j \in H^1(0,l_j) \ \forall j \in J$  and $y|_{\Gamma}=0.$ 
		Let $\mathcal{H}^{-1}$ be the dual space of $\mathcal{H}^1$.
		We further introduce the space $ \mathcal{H}^2$ of  functions $y \in \mathcal{H}^1$  such that }
		$y_j \in H^2(e_j)$ for every $j \in J$ and,
		for every $v_i \in V \setminus \Ga$, the equality
		\begin{equation}\label{s-d}
			\sum_{j \in J(v_i)} \partial y_j(v_i)=0
		\end{equation}
		holds. Here (and everywhere later on)
		$\partial y_j(v_i)$ denotes the
		derivative of $y$ at the vertex $v_i$ taken along the  edge  $e_{j}$
		in the direction outwards the vertex.
		Vertex conditions \eqref{s-d} (together with continuity at $v$) are called the 
		 {\bf Kirchhoff-Neumann} (KN) conditions for the internal vertices. 
		
		Let $q$ be a real valued function (potential) such that $q|_{e_j} \in C[0,l_j]$. We define the  self adjoint (Schr\"odinger) operator $A$ on the graph $\Om$  
		whose operation is 
		\begin{equation} \label{sch}
			(Ay)_j  = - y_j''(x) + q_j(x)y_j(x), \ \ \forall j \in J, \,
		\end{equation}
		in $ \mathcal{H} $ with the operator domain $\mH^2$.

	{Consider  the following initial boundary value problem (IBVP): 
		
		\begin{equation} \label{waveo}
			w_{tt}-w_{xx} +q(x)w=p(x,t) \quad \textrm{ in} \ \  \{\Omega\setminus V\} \times
			(0,T),
		\end{equation}
		\begin{equation} \label{kco1}
			w_j(v_i,t)=w_k(v_i,t) \   {\rm for}  \ j,k \in J(v_i), v_i \in V\setminus \Gamma,  \ t \in [0,T],
		\end{equation}
		\begin{equation}\label{kco2}
			\sum_{{j} \in J(v_i)} \partial w_j(v_i,t)=0  \ \ {\rm at \ each \ vertex } \ \ v_i \in V\setminus \Gamma, \ t \in [0,T],\\
		\end{equation}
		\begin{equation}\label{kco3}
			w_j(v_i,t)=0  \ \ {\rm at \ each \ vertex } \ \ v_i \in \Gamma, \ t \in [0,T],\\
		\end{equation}
		\begin{equation} \label{ico}
			w|_{t=T}=w^0, \ w_t|_{t=T}=w^1 \ \ \textrm{ in} \ \ \Omega.
		\end{equation}
		Here $T>0, \, p \in L^1(0,T;\mathcal{H}), \, w^0 \in \mathcal{H}^1,$ and  $w^1 \in \mathcal{H}.$    Using the Fourier method, one can show (see Lemma \ref{dual} below) that this IBVP has a unique generalized solution such that for any $i \in I, j \in J$,
		
		\begin{equation} \label{reo}
			w \in C([0,T];\mathcal{H}^1), \ w_t \in C([0,T];\mathcal{H}), \ w(v_i,\cdot)\in H^1(0,T),\ \partial w_j(v_i,\cdot) \in L^2(0,T) 
		\end{equation}

		We now parametrize a set of observations, using sets $I^*\subset I\setminus \Gamma$,  and $J^*\subset J$.
We denote $V^*=\{ v_i\in V,\ i\in I^*\}.$
		If $j\in J^*$,
		we define $i=\iota (j)$ to be one of the vertices at the end of $e_j$;\je{ here $\iota$ is a function from $J^*$ to $I$.} Define
		$J^*(v_i)$ to be the set $ J^*\cap J(v_i)$, and $J^c(v_i)=J(v_i)\setminus J^*.$ 
  We will call $\{ I^*,J^*,\iota\}$ the active set. Our observations will be 
  $
  w(v_i,t),\ i\in I^*,$ and $\pa w_j(v_i,t)$, with $j\in J^*(v_i),\ i=\iota (j).$
 The observability inequality is then 
  $$\| w^0\|_{\mathcal{H}^1}^2+\| w^1\|_{\mathcal{H}}^2
 +\| p\|^2_{L^1(0,T;\mathcal{H}^1)} 
  \leq C (\sum_{i\in I^*}\int_0^T|w(v_i,T-t)|^2dt+
  \sum_{j\in J^*,i=\iota (j)}\int_0^T|\pa w_j(v_i,T-t)|^2dt).
  $$
  Here, $w^0,w^1$ are arbitrary, and $C$ is a constant independent of $w^0,w^1$.  We emphasize that our observations are feasible and nondestructive,
and hence the corresponding exact control problem will be natural.

We can now define the system we wish to control.
	The control system dual to \eqref{waveo}--\eqref{ico} with the active set $\{I^*, J^*,\iota \}$ can be written as follows:
		\begin{equation} \label{wave1}
			u_{tt}-u_{xx} +q(x)u=0 \quad \textrm{ in} \ \  \{\Omega\setminus V\} \times
			(0,T),
		\end{equation}
		
		\begin{equation} \label{nc}
			\sum_{j \in J(v_i)}  \partial u_j(v_i,t)= 
			\begin{cases} 
				0, & v_i\notin \Gamma, \ i \in I \setminus I^{*}, \\ \tilde{f}_i(t), & v_i\notin \Gamma,\  i \in I^{*}, \,  \end{cases}
		\end{equation}
		
		\begin{equation} \label{bc}
			u_j(v_i,t) = f_i(t), \je{\ v_i\in \Gamma, \ i=\iota (j),\ j\in J^*,} 
		\end{equation}
		\begin{equation} \label{bc2}
			u_j(v_i,t) = 0, \ \je{\ v_i\in \Gamma, \ i=\iota (j),\ j\notin J^*, }
		\end{equation}
		\begin{equation} \label{dc}
			\begin{cases} 
				u_j(v_i,t)-u_k(v_i,t)=f_j(t), & i \in I\setminus \Gamma, \; j \in J^*(v_i), \; k \in J^c(v_i), \\
				u_j(v_i,t)-u_k(v_i,t)=0, & i \in I\setminus \Gamma, \; j, k \in J^c(v_i),  \end{cases}
		\end{equation}
		
		\begin{equation} \label{ic1}
			u|_{t=0}=u_t|_{t=0}=0 \ \ \textrm{ in} \ \ \Omega.
		\end{equation}
		In this paper we assume that $f_i \in L^2(0,T) $, and 
		$\tilde{f}_i\in H^1(0,T)'$. 
		We will refer to $f_i(t)$ with $i\in I^*$, 
as a ``vertex control", and to $f_j(t)$ with $j\in J^*$ as an ``edge control". We will state our controllability result in the next section.

		\ 
		
	We conclude this section by defining the relevant spectral data.
		Let $\{ (\omega^2_n,\f_n): n\geq 1\}$ be the eigenvalues and normalized eigenfunctions 
		of the operator $L$, i.e. $\f_n$ solves the following eigenvalue problem:
		
		\beq \lab{sp}
		- \f_n'' + q(x) \f_n = \om_n^2 \f_n \quad \rm{on} \ \ \{\Omega\setminus V\},
		\eeq
		\beq \lab{ncf} 
		\sum_{j \in J(v_i)} \pa (\f_n)_j(v_i)=0, \  v_i \in V \setminus \Ga, \ \f_n|_{\Ga}=0,
		\eeq
		\beq \lab{kn}
		(\f_n)_j(v_i)=(\f_n)_k(v_i), \quad  j, k \in J(v_i),\ v_i \in V \setminus \Ga.
		\eeq
		where $(\f_n)_j$ is the restriction of $\f_n$ to edge $e_j$. 
		
		\begin{subsection}{Well-posedness, Part 1}
			The main goal of this section is to determine the spaces for which the system \rq{wave1}-\rq{ic1} is well-posed.
			{For general graphs  with $L^2$ Dirichlet controls applied at the boundary, it is well known (see, e.g. \cite{LLS,AN}) } that 
			the mapping $t\mapsto u(\cdot,t)$ is in 
			$C(0,T;\mH)\cap C^1(0,T;\mH^{-1}).$
			However, this is not the case for controls applied in the interior, as our next example shows. 
			
			\begin{example}
				Consider the star graph consisting of three copies of $(0,l)$ with central vertex identified with $x=0$ for each edge. We set $q=0$, and assume a single Neumann control 
				$\tilde{f}(t)$ is applied at $x=0$. Letting $g(t)=u(0,t)$, for $t<l$ we have $u_j(x,t)=g(t-x)$ so that 
				$$\tilde{f}(t)=\sum_{j=1}^3\pa_xu_j(0,t)=3g'(t).
				$$
				Let $\phi\in \mH^1_0$. Then assuming for the moment that $\tilde{f}$ is sufficiently regular, from integration by parts we have 
				\begin{eqnarray*}
					\int_{\Omega}\phi u_t & = & \sum_j\int_0^l\phi_j(x)g'(t-x)dx\\
					& = & \sum_j  \int_0^l\phi_j'(x)g(t-x)dx +3\phi (0)g(t).
				\end{eqnarray*}
    {I propose to use $f(\cdot),$ not  $f(*),$ below and everywhere}
				Thus $\< \phi (\cdot) ,u_t(\cdot,t)\>_{1,-1}$ being continuous in $t$ for $g\in L^2(0,T)$ requires $\phi (0)=0.$ {Here $\< \cdot,\cdot\>_{1,{-1}}$ is the pairing $\mH^1-\mH^{-1}.$}
			\end{example}
			With this example in mind, we parametrize the internal vertices associated either to an edge or vertex control. 
			{We define $\bar{I}=\{ i\in I: i\in I^* \mbox{ or } J^*(v_i)\neq \emptyset \} $
				and} $$
			\mH^1_b=\{ \phi\in \mH^1_0: \ \phi(v_i)=0 \ \forall \, i\in \bar{I} \}, \ \  \mH^2_{b}=\{ \phi\in \mH^1_b: \phi |_{e_j}\in H^2(e_j), \ \forall \, j \in J
			\},$$
		$$\mH^1_0=\{ \phi\in \mH^1_0: \ \phi(v_i)=0 \ \forall \, i\in {I}^* \}. $$
			{In this section, we will prove that the velocity belongs to $(\mH^1_0)',$  the dual space of $\mH^1_0.$ } We remark in passing that for such a velocity, the Dirac delta function $\delta_{v_i}$ would be non-trivial only if $i \in I \setminus {I}^*. $
			
			\je{We also define
			$$H^1_0 (0,T)=\{ g\in H^1(0,T): \ g(0)=g(T)=0\}, \ \ H^2_0 (0,T)=\{ g\in H^2(0,T): \ g(0)=g'(0)=g(T)=g'(T)=0\}.
			$$
			We will use two kinds of controls: ``smooth'' and ``rough''.
				For all edge controls, including boundary, $f_j\in H^2_0 (0,T)$ for smooth controls and $f_j\in L^2 (0,T)$ for rough controls, while all
				vertex controls will satisfy 
				$\tf_i\in H^1_0(0,T)$ or $H^1(0,T)',$ correspondingly.
			
			The corresponding control spaces for the system \rq{wave1}-\rq{ic1} are: 
			$$
			\F_2^T=(H_0^1(0,T))^{|I_*|}\times (H_0^2 (0,T))^{|J^*|};  \ \
			\F_0^T=(H^1(0,T)')^{|I_*|}\times (L^2 (0,T))^{|J^*|}.
			$$}

			\begin{thm}\label{thm1}
				\ 
				
				A) Given ${\bf f}\in \F_2^T$, there exists a unique solution $u^{\bf f}$ solving \rq{wave1}-\rq{ic1}, such that $$u^{\bf f}(x,t)\in C(0,T;\mH^2)\cap C^1(0,T:\mH^1).
				$$
				
				B) Given ${\bf f}\in \F_0^T$, there exists a unique solution $u^{\bf f}$ solving \rq{wave1}-\rq{ic1}, such that $$u^{\bf f}(x,t)\in C(0,T;\mH)\cap C^1(0,T:(\mH_0^{1})').
				$$
			\end{thm}
		The proof of part A, where the solution $u$ will be classical, is left to the reader.

		\je{We begin the proof of  part B by proving the following, slightly weaker result.
	\begin{prop}	\label{half}
	Given ${\bf f}\in \F_0^T$, there exists a unique solution $u^{\bf f}$ solving \rq{wave1}-\rq{ic1}, such that $$u^{\bf f}(x,t)\in C(0,T;\mH)\cap C^1(0,T:(\mH_b^{1})').
				$$	
	\end{prop}	}
Proof of proposition: 		

We consider the weak formulation of the problem. Thus	
			we now consider the problem  \rq{waveo}-\rq{ico} with 
			$w_0=w_1=0$.
			In what follows, when convenient, we will denote the function $w$ determined { by $p$ by $w^p.$}
			Assuming classical solutions, multiplying \rq{waveo} by $u$ and then integrating by parts,
			we get
			\begin{eqnarray}
				\int_Qpu\ dxdt & =   & \int_0^T \left[\sum_{ j \in J^*(v_i), i=\iota (j)} \partial w_j(v_i,T-t) \, f_j(t)
				-\sum_{i \in I^*} w(v_i,T-t) \,  \tilde{f}_i(t)  \right] \, dt .\label{ide0}
			\end{eqnarray}

			Thus, for a given ${\bf f}\in \F^T_0$,  we will define $u$ to be the weak solution as the solution of  the system \rq{wave1}-\rq{ic1} if
			\begin{eqnarray}
				\int_Q p(x,t)u(x,t)\ dxdt & =   & \int_0^T \sum_{ j \in J^*, i=\iota (j)} \partial w_j(v_i,t) \, f_j(t)\ dt \ 
				-\sum_{i \in I^*}\<  w(v_i,t) \,  \tilde{f}_i(t)\>_{1,-1}  .\label{weakwg}
			\end{eqnarray}
			holds for all $p\in { L^1(0,T; \mH)}$. Here $\< *,* \>_{1,-1}$ is the pairing $H^1(0,T)-H^1(0,T)'$.
			
			We now consider the regularity of $w^p$.
			\begin{lemma}\label{dual}
				{ Assume that $q_j \in C^1[0,l_j]$ and $q$ is continuous at the interior vertices.}

				A) Suppose $p\in L^1(0,T; \mH)$. Then
				$w^p\in C(0,T; \mH^1_0)\cap C^1(0,T;\mH).$
				Furthermore, for all $v_i\in V$,
				\beq
				\int_0^T 
				| w_j(v_i,t)|^2dt \leq C \|p\|^2_{ L^1(0,T; \mH )}.
				\label{e1}
				\eeq
				Finally, for all \je{$i\in I,$} \je{and all $j\in J(v_i),$}
				\beq
				\int_0^T|\je{\pa w_j(v_i,t)|^2}+|w_t(v_i,t)|^2 \leq C \|p\|^2_{ L^1(0,T; \mH )}.
				\label{e2}
				\eeq
				Here $C$ is independent of $p$. 
				
				B) Suppose $p=\frac{\pa}{\pa t}P$, with 
				$P\in L^1(0,T; \mH^1_b)$.
				Then
				$w^p\in C(0,T; \mH )\cap C^1(0,T;(\mH_b^{1})').$
				Furthermore, for all $v_i\in V$,
				\beq 
				\int_0^T 
				|\je{w}(v_i,t)|^2dt \leq C \|P\|^2_{ L^1(0,T; \mH_b^1 )}.\label{e3}
				\eeq
				Finally, for all \je{$i\in \overline{I},$ and all $j\in J(v_i)$}
				\beq
				\int_0^T|\je{\pa w_j}(v_i,t)|^2+|w_t(v_i,t)|^2 \leq C \|P\|^2_{ L^1(0,T; \mH_b^1 )}.\label{e4}
				\eeq
				Here $C$ is independent of $f$.
			\end{lemma}
			Proof: Let $Q=\Omega \times [0,T].$
			$\| w_t(*,t)\|_{L^2(\Omega)}+\| w_x(*,t)\|_{L^2(\Omega)}$
			\je{To prove the continuity results and \rq{e1}, we use a Fourier series argument. 
			Recall we denote the eigenvalues and unit eigenfunctions for our system by $\{ (\om_j^2,\f_j\}_1^{\infty}.$ Here the eigenvalues are listed in increasing order. Assume for the moment that $\la_j^2>0$ for all $j.$
			Letting $p=\sum_jp_j(t)\f_n(x)$, standard calculations give
		\beq
			w(x,t)=\sum_{j=1}^\infty \big (\frac{1}{\om_j}\int_0^t\sin (\om_j(t-s))  \ ds \big )\f_j(x).
		\label{w}
		\eeq
			In what follows, assume first that $p$ is regular, and denote $\| *\|$ the $L^2$-norm on $\Omega$. Then 
			\begin{eqnarray*}
			\| w_x(*,t)\|^2 & = & \< \pa_x^2w(*,t),w(*,t)\>_{{\cal H}}\\
			& - &\sum_j  (\int_0^t\sin (\om_j(t-s))p_j(s)  \ ds)^2\\
			& \leq & {{t}}\big ( \sum_j   \int_0^t|p_j(s)|^2  \ ds\big )\\
			& \leq & {{t}}    \int_0^t\sum_j|p_j(s)|^2  \ ds.
			\end{eqnarray*}
			Hence,
			\begin{eqnarray*}
			\| w_x(*,t)\|	& \leq & {\sqrt{T}}  \int_0^T\big (  \sum_j|p_j(s)|^2  \big )^{1/2}\ ds\\
				& \leq & {\sqrt{T}} \| p\|_{L^1(0,T;{\cal H})},\ t\in [0,T].
			\end{eqnarray*}
			We conclude 
			$$\| w_x\|_{L^{\infty}(0,T;{\cal H})}\leq C  \| p\|_{L^1(0,T;{\cal H})}.$$
			It follows by a density argument that $w\in C(0,T;\mH^1_0)$. A similar argument shows 
			$w\in C^1(0,T;\mH)$. 	This proves \rq{e1}.
			
			The argument above can easily be adapted to
			 the case where for some $j$ we have $\om_j^2\leq 0.$  For $\om_j<0$, it suffices to replace $\sin (\om_j(t-s))$ by $\sinh (|\om_j|(t-s))$, and for $\om_j=0$ we can replace 
			$\sin (\om_j(t-s))$ by $(t-s)$.}

			To prove \rq{e2}, we use a multiplier argument. Fix any vertex $v_i$.
			Now parametrize the incident edges by $(0,l_j)$ with $v_i$ corresponding to $x=0.$  Let $h(x)$ be some smooth cut-off function with $h(0)=1,$ $h_j'(0^+)=0,$ and $h(x)=0$ for \je{$x >\min (l_j/2),$} and set $h=0$ elsewhere on $\Omega,$ and further assume $h'(x)\leq 0.$ Then
			\begin{eqnarray*}
				\int_Qhw_xp& = & \sum_{j\in J(v_i)}\int_0^T\int_{e_j}
				hw_x(w_{tt}-w_{xx}+qw)\ dxdt\\
				& = & \frac{-1}{2}\sum_{j\in J(v_i)}\int_0^T\int_0^{l_j}
			h\pa_x(w_{t}^2+w_{x}^2)-hq\pa_x(w^2)\ dxdt\je{+\int_{\Omega}[hw_xw_t]_0^Tdx}\\
				& = & \frac{1}{2}\sum_{j\in J(v_i)}\int_0^T\int_0^{l_j}
				h_x(w_{t}^2+w_{x}^2)-(hq)_xw^2\ dxdt\je{+\int_{\Omega}[hw_xw_t]_0^Tdx}
				\\
				& & +\frac{1}{2}\sum_{j\in J(v_i)}\int_0^Tw_t(0,t)^2+\pa_xw_j(0,t)^2-q_j(0^+)w(0,t)^2 \, dt.
			\end{eqnarray*}
			\je{The estimate \rq{e2} now follows from $w\in C(0,T;\mH^1_0)\cap C^1(0,T;\mH)$.}

			\je{ We now prove part B of the Lemma, which is an adaptation of Lions (\cite{Li}, Thm. 4.2, p. 46).} Let $\theta$ solve
			\begin{eqnarray}
				\theta_{tt}-\theta_{xx} +q(x)\theta & =&P(x,t) \quad \textrm{ in} \ \  \{\Omega\setminus V\} \times
				(0,T),\label{theta} \\
				\theta|_{t=T} =  \ \theta_t|_{t=T} &  = & 0,  \\
				\theta_j(v_i,t)& =& \theta_k(v_i,t) \   {\rm for}  \ j,k \in J(v_i), v_i \in V\setminus \Gamma,  \ t \in [0,T],\label{kn1} \\
				\sum_{{j} \in J(v_i)} \partial \theta_j(v_i,t)&=&0  \ \ {\rm at \ each \ vertex } \ \ v_i \in V\setminus \Gamma, \ t \in [0,T],\label{kn2} \\
				\theta_j(v_i,t)& = &0  \ \ {\rm at \ each \ vertex } \ \ v_i \in \Gamma, \ t \in [0,T].
			\end{eqnarray}
			{ The inclusion } $\theta\in C(0,T;\mH)$ can be proven using the multiplier $h \theta_x$; the details are left to the reader.
			
			Clearly $\theta_t=w.$ Let $v_i$ be any interior vertex.
			\je{We wish to prove an upper bound for} 
			\beq
			\int_0^Tw(v_i,t)^2+w_t(v_i,t)^2+\pa_xw_j(v_i,t)^2\ dt=\int_0^T\theta_t(v_i,t)^2+\theta_{tt}(v_i,t)^2+\pa_x
			\theta_t(v_i,t)^2\ dt. \label{kest}
			\eeq 
				First, 		using \rq{kn1},\rq{kn2},  one can check that for sufficiently smooth $\theta,$
				$$
				\int_\Omega \theta_{tt}\theta_{xxt}=-\int_\Omega \theta_{ttx}\theta_{xt}.
				$$
			Assume for the moment that $P$ is regular in $t$ and vanishes near $t=0,T.$
			We multiply \rq{theta} by $\theta_{xxt}$ and integrate by parts as above over $\Omega$ to get
			\begin{eqnarray*}
				\frac{d}{dt}\int_\Omega (\theta_{xx}^2+\theta_{xt}^2 \ dx)=\int_\Omega  (P_x\theta_{xt}\je{-}q_x\theta\theta_{xt}\je{-}q\theta_x
				\theta_{xt})\, dx+\je{\sum_{v_i\in V\setminus \Gamma}
			\left (	\big ( P(v_i,t)-q(v_i)\theta(v_i,t)\big )
				\sum_{j\in J(v_i)}(\theta_j)_{xt}(v_i,t)\right )}.
			\end{eqnarray*}
			The boundary terms on the right hand side of the equation above vanish because of the continuity of $q$ at interior vertices \je{and \rq{kco2}.}
			A simple exercise in analysis now gives
			\beq
			\| \theta\|_{C(0,T;\mH^2_{b} )\cap C^1(0,T;\mH^1_{b})}\leq C \| P\|_{L^1(0,T; \mH^1_{b})}.\label{cont}
			\eeq
			By density, this estimate will hold for all $P\in L^1(0,T; \mH^1_{b}).$
			\je{Equation \rq{e3} now follows.}
			
	\je{	We now prove \rq{e4}, so assume $i\in \overline{I}.$}	
			For the moment in this paragraph, we assume 
			$P$ is regular in $t$ and vanishes at $t=0,T.$
			In what follows we will use $h$ as in the proof of Part A to localize to the edges incident to $v_i$.
			We  multiply \rq{theta} by $\theta_{xtt}(x,t)h(x)$ and integrate over  $Q$ to get
			\begin{eqnarray}
				\int_Q \theta_{ttx}hP \ dxdt & = &
				-\frac{1}{2}\int_Qh_x(\theta_{tt}^2+\theta_{tx}^2)\ dxdt-\frac{1}{2}\int_0^T\theta_{tt}(v_i,t)^2+\theta_{xt}(v_i,t)^2dt\nonumber\\
				&  &+
				\int_\Omega\theta_{tx}(x,0)\theta_{xx}(x,0)h(x)\
				dx-\int_Q\theta_{xt}hq\theta_t\ dxdt\nonumber\\
				& & -\int_\Omega\theta_{xt}(x,0)h(x)q(x)\theta(x,0)\ dx.
				\label{multi}
			\end{eqnarray}
			Integration by parts on the left hand side, together with 
			$\theta_{tt}=\theta_{xx}-q\theta+P$, gives
			$$\int_Q \theta_{ttx}hP \ dxdt=-\int_Q\theta_{xx}(hP)_x-q\theta(hP)_x +\frac{h_x}{2}P^2\ dxdt+\sum_{j\in J(v_i)}\int_0^T (-\theta_{tt}(v_i,t)P(v_i,t)+\frac{1}{2}P(v_i,t)^2)\ dt .$$
			The last two terms vanish because  $P(v_i,t)=0$, $i\in \overline{I}.$
			\je{Remark: it is precisely in the last equation where the assumption $P(v_i,t)=0$  seems essential. }
			Inserting this along with $\theta_{tt}=\theta_{xx}-q\theta+P$ into \rq{multi} gives
			$$
			\frac{1}{2}\int_0^T\theta_{tt}(v_i,t)^2+\theta_{xt}(v_i,t)^2\ dt = -\int_Qh_x\big( \theta_{xx}(-q\theta +P)+\frac{1}{2}(q^2\theta^2+\theta_{xx}^2+\theta_{tx}^2)\ dxdt+\int_\Omega\theta_{tx}(x,0)\theta_{xx}(x,0)h(x)\ dx
			$$
			$$-\int_Q(q\theta(Ph)_x+\theta_{xt}hq\theta_t-\theta_{xx}(hP)_x)\ dxdt
			-\int_\Omega\theta_{xt}(x,0)hq(x)\theta(x,0)\ dx.
			$$
			Thus by \rq{cont}, 
			$$\int_0^T\je{\theta_{tt}(v_i,t)^2+\theta_{tx}(v_i,t)^2dt}=
			\int_0^Tw_{t}(v_i,t)^2+w_{x}(v_i,t)^2 dt\leq 
			\| P\|_{L^1(0,T; \mH^1_b)}.
			$$
			This inequality was proven for $P$ regular in $t$ and vanishing at $t=0$, but by density it will hold for $P\in {L^1(0,T; \mH^1_{b})}. $
			The lemma is proven.
			
			\

		\je	{\bf Proof of Proposition \ref{half}}.
			
		\je{We write 
			\beq
			{\bf f}=
			\< \tilde{f}_1,...,\tilde{f}_{|I^*|},f_1,...,f_{|J^*|}\>.
			\label{g}
			\eeq
			By linearity, it suffices to prove \rq{weakwg} when all but one component is zero. We prove the result for ${\bf f}=\< \tilde{f}_1,0,...,0\>$, leaving the result for Dirichlet controls to the reader. 
			For readability, we rewrite $\< \tilde{f}_1,0,...,0\>$ as $g(t)\in H^1_*(0,T),$ and assume the vertex associated to $g$ is $v_1$.	}
			
			By Part A of Lemma \ref{dual}, the left hand side of \rq{weakwg} is a bounded linear functional on $L^1(0,T; \mH )$. Thus,
			there exists a unique $u\in L^{\infty}(0,T; \mH )$ such that \rq{weakwg} holds, and furthermore there exists $C$ such that 
			\beq
			\sup_t \|u(*,t)\|_{\mH }\leq C\|g\|_{H^1(0,T)'}.\label{est1}
			\eeq
			That $u \in C(0,T; \mH )$ follows by the following argument
			density argument. Clearly for  regular $g_n$, $u^{g_n}\in C(0,T;\mH ).$ Suppose $g_n\to g$ in $H^1(0,T)'.$
			It follows by \rq{est1} that 
			$$\sup_t \| u^g(x,t)-u^{g_n}(x,t)\|_{\mH }\to 0.
			$$
			Continuity of $t\mapsto u^g$ follows. 
			
			We now prove continuity of $t\mapsto u_t^g$.
			Assume for the moment $P\in C_0^{\infty}(0,T; \mH^1_b)$ and hence vanishes at $t=0,T$. Setting $p=\pa_tP$ in \rq{weakwg} and then integrating  by parts, we get
			$$
			\<g(*), w^p(v_1,*)\>_{-1,1}=-\int_0^T\int_\Omega Pu_t\ dxdt.
			$$
			By Lemma \ref{dual} part B, the left hand side of this equation defines a bounded linear functional on $L^1(0,T;\mH^1_b)$, and hence $u_t\in L^{\infty}(0,T; (\mH^1_b)')$ and 
			$$
			\sup_t \|u_t(*,t)\|_{(\mH^1_b)'}\leq C\|g\|_{H^1(0,T)'}.
			$$
			Continuity follows from another density argument. 
			
			\

			\ 
			
			\je{The following technical result will be important in proving our exact control result.}
			\begin{prop}
				Let ${\bf f}\in {\cal F}^T_2.$ Then 
				\beq
				\pa_t^2 u_t^{\bf f}(*,t)=u_t^{{\bf f}''}(*,t), \label{ident}
				\eeq
				this equation holding in $C(0,T;(\mH_0^1)').$
			\end{prop}
			We remark this equation holds trivially for ${\bf f}$
			sufficiently regular. Also, for component $\tilde{f}_i\in {H}^1_0(0,T),$ 
			\je{we interpret $\tilde{f}_i''$ as follows. Defining $D_t: H^1(0,T)\mapsto L^2(0,T)$ by $D_th(t)=h'(t)$, we have $D_t^*: L^(0,T)\mapsto H^1(0,T)'$. Then  $$\tilde{f}_i''=-D_t^*D_t\tilde{f}_i\in H^1(0,T)'.$$}
			
			Proof of proposition: 
			\je{We define $g(t)$ as in the proof of Theorem \ref{thm1}, see \rq{g}. For this proof, we denote $C_*^1(0,T)=\{ g\in C^1(0,T): \ g(0)=0\}.$}
			It suffices to prove
			$$\pa_t^2 u^g(x,t)=u^{g''}(x,t)
			$$
			as elements of $C(0,T;\mH ).$ 
			Since $g''\in H^1(0,T)'$, we have by Theorem 1 that 
			$u^{g''}$ solves \rq{wave1}-\rq{ic1}.
			We now prove 
			$$
			\< w^p(v_1,*),g''(*)\>_{H^1(0,T),(H^1)'(0,T)} =
			-\int_Q \pa_t^2 u^{g}(x,t)p(x,t) \ \je{dxdt},\ \forall p\in { L^1(0,T; \mH )};
			$$
			the proposition will then follow by \rq{weakwg} and uniqueness of the weak solution. 
			Assume for the moment $p\in C_0^{\infty}(0,T; \mH^1_0)$, and hence
			$$w^p_{xx}(x,T)=w^p_{xxt}(x,T)=0,
			$$
			and $w^p_{ttxx}\in C(0,T;\mH ).$ Also assume for the moment that $g\in C^1_*(0,T)$, so that $u^g$ has a classical solution. In the argument that follows, we will integrate by parts in $x$.
			Thus, we assume that $\Omega$ is a directed graph, and for each edge we have the identification
			$e_j=\{ x\in (a_j,b_j)\}$, with $b_j$ always chosen to be the incoming vertex.

			Hence, 
			\begin{eqnarray}
				\int_Q\pa_t^2 u^{g}(x,t)p(x,t) \ dxdt& = & 
				\int_Q u^{g}_{tt}(w_{tt}-w_{xx}+qw) \ dxdt\nonumber\\
				& = & 
				\int_Q  u^{g}_{tt}w_{tt}-u^gw_{xxtt}+qu^gw_{tt} \ dxdt\nonumber\\
				& = & 
				\int_0^T\sum_{j\in J}[ -u^{g}w_{ttx}+u^g_xw_{tt}]_{a_j}^{b_j} \ dt\nonumber\\
				& = & 
				\int_0^T\sum_{i\in I}\sum_{j\in J(v_i)}
				\big( (u^{g}_j(v_i,t)(w_j)_{ttx}(v_i,t)
				-(u^g_j)_x(v_i,t)(w_j)_{tt}(v_i,t)\big) \ dt\nonumber\\
				& = &-
				\int_0^Tg(t)w_{tt}(v_1,t)\ dt\nonumber\\
				& = &
				\int_0^Tg'(t)w_{t}(v_1,t)\ dt\nonumber\\
				& = &
				\<D_tg(*),D_tw(v_1,*)\>_{L^2(0,T),L^2(0,T)}\nonumber\\
				& = & -\<g'',w(v_1,*)\>_{H^1(0,T)',H^1(0,T)}.\label{id}
			\end{eqnarray}
			This identity is proven for all $p\in C_0^{\infty}(0,T; \mH^1_0) $ and $g\in C_*^1(0,T)$. We now use two density arguments. An exercise in Fourier series shows that for $g\in \je{{H}^1_*(0,T)},$ we have 
			$$
			\| \pa^2_tu^g\|_{L^2(0,T;\mH )}\leq C \| g\|_{H^1(0,T)}.
			$$
			Thus letting $g_n\in C_*^1(0,T)$ approach $g$ in $H^1(0,T)$, we see that \rq{id} holds for $g\in \je{{H}^1_*(0,T)}.$
			Now
			let $p_n \in C_0^{1}(0,T; \mH) $ be such that $p_n\to p$ in \je{$L^1(0,T; \mH^1_0)$}. Thus
			we have $w^{p_n}(v_1,t)\to w^p(v_1,t)$ in $H^1(0,T) $ by Lemma \ref{dual}, part A.  Since $\pa^2_tu^g\in L^\infty (0,T; \mH)$, we conclude \rq{id} holds for all $p$ in $L^1(0,T; \mH )$. The proposition is proved.
		\end{subsection} 
		
	\subsection{The forward problem on an interval}  \label{sec:edge}

In this section we discuss two IBVP problems of the wave equation on an interval $[0,l]$ with  Dirichlet boundary conditions. The solutions will be used in solving  forward and control problems on general graphs later.  We first consider the following problem:

\begin{equation} \label{stdwaveint}
\begin{cases} (u)_{tt}-(u)_{xx}+q(x)u=0, \quad  0<x<l, \ 0<t<T\\u|_{t \le 0} = 0, \quad  u(0,t)=f(t), \quad u(l,t)=0. \end{cases}
\end{equation}
We will refer to $f$ as the Dirichlet control function or simply, Dirichlet control.  For $f \in H_*^1(0,T)$, the system \eqref{stdwaveint} has a unique solution $u^{f, DD} \in C([0,T]; H^1(0,l)).$ 
This solution can be presented by the so-called  folding ruler formula \cite{AZ}:
\begin{eqnarray}
u^{f,DD,0}(x,t)  &
=  &  f(t-x) + \int_x^t \omega(x,s)f(t-s) \, ds \nonumber\\
&& -f(t-2l+x)-\int_{2l-x}^t \omega(2l-x,s)f(t-s) \, ds \nonumber \\
&&+f(t-2l-x)+\int_{2l+x}^t \omega(2l+x,s)f(t-s) \, ds \nonumber \\
&&-f(t-4l+x)-\int_{4l-x}^t \omega(4l-x,s)f(t-s) \, ds + \dots   \nonumber \\
&= &  \sum_{n=0}^{\left \lfloor{\frac{t-x}{2l}}\right \rfloor} \left( f(t-2nl-x) +\int_{2nl+x}^t \omega( 2nl+x,s) f(t-s) \, ds \right) \nonumber \\
&& -\sum_{n=1}^{\left \lfloor{\frac{t+x}{2l}}\right \rfloor} \left( f(t-2nl+x) + \int_{2nl-x}^t \omega( 2nl-x,s) f(t-s) \, ds \right)  \label{fold2} 
\end{eqnarray}
\noindent where  $\left \lfloor{\cdot}\right \rfloor$ is the floor function;  $\omega(x,t)$ is a solution to the Goursat problem 
\begin{equation} \label{omega}
\begin{cases} (\omega)_{tt}-(\omega)_{xx}+q(x)\omega=0, \quad  0<x<\infty, \\
 \omega(0,t)=0, \omega(x,x)=-\frac{1}{2}\int_0^x q(s) \, ds, \end{cases}
\end{equation}
in which the potential $q(x)$ is extended to the semi-axis $x>0$ by the rule $q(2nl \pm x)=q(x)$ for all $n \in \mathbb{N}$.  

The folding ruler formula gives convenient presentations of the solution with various boundary conditions and will be used in many constructions of this paper. More details about this formula can be found in \cite{AZ}. 
When the Dirichlet control function $f \in H_*^1(0,T)$ is applied at $x=l$, the IBVP 
\begin{equation} \label{stdwaveintbk}
\begin{cases} (u)_{tt}-(u)_{xx}+q(x)u=0, \quad  0<x<l, t>0 \\ u|_{t \le 0} = 0, \quad  u(0,t)=0, \quad u(l,t)=g(t). \end{cases}
\end{equation}
can be solved by changing of variables in \eqref{stdwaveint}.  We put $p(x)=q(l-x)$ and extend $p$ by letting $p(2nl \pm x)=p(x)$. Let $\bar{\omega}(x,t)$ be the solution to the Goursat problem \eqref{omega} where $q(x)$ is replaced with $p(x)$, then 
\begin{align} \label{fold2bk} 
u^{0,DD,f}(x,t) \quad  = \quad &f(t-l+x) + \int_{l-x}^t \bar{\omega}(l-x,s) f(t-s) \, ds   \nonumber \\
&- f(t-l-x) -\int_{l+x}^t \bar{\omega}(l+x,s) f(t-s) \, ds +\dots \nonumber \\
\end{align}

When Dirichlet control functions $f(t)$ and $g(t)$ from $H^1_*(0,T)$ are applied at $x=0$ and $x=l$ respectively, the solution to for the wave equation \eqref{wave1} on $e$ with zero initial condition is denoted as $u^{f,DD,g}$.  By the superposition principle, $u^{f,DD,g}(t) = u^{f,DD,0}(t)+u^{0,DD,g}(t)$.  

\subsection{Well-posedness, Part 2}
\je{					In this section, we complete the proof of Theorem \ref{thm1}
		
We need to sharpen our results about the regularity of $u_t$
at interval vertices $v_i$ for which no Neumann control is applied, i.e. $i\notin I^*.$  By locality, 
it suffices to consider the behaviour in a small neighborhood of $v_i$, so in what follows
we can consider a star graph where the edges $\{ e_j\}_1^N $ are identified with the intervals $\{ (0,l_j)\}$  and the
central vertex $v_i$ is identified with $x=0$. By linearity it is sufficient to assume there is a single Dirichlet control at $v_i$.
Let $u$ be the solution to the system \rq{wave1}-\rq{ic1}, and set $g_j(t)=u_j(0,t)$.

Assume the single control $f$ is applied at $x=0$ on edge $e_2$.
Thus
\beq \label{vc}
g_2(t)=g_1(t)+f(t), \ g_j(t)=g_1(t), \ j=3,...,N.
\eeq 
For the moment, assume $f$, and hence $g_j$, is regular.
Suppose for the moment that $T$ is small, so for $t<T$ we can disregard any reflected waves, and 
$$u_j(x,t)=g_j(t-x)+\int_x^tk_j(x,s)g_j(t-s)ds.$$
Thus the equation $\sum \pa u_j(0,t)=0$ implies
$$0=\sum_j g_j'(t)-\int_0^t \pa k_j(0,s)g_j(t-s)ds.
$$			
Integrating this, and using $g_j(0)=0$, we get
\beq \label{KN'}
0=\sum_j g_j(t)-\int_0^t \int_0^r\pa k_j(0,s)g_j(r-s)\ dsdr.
\eeq
Hence 
$$
\sum_j \big ( g_1(t)-\int_0^t \int_0^r\pa k_j(0,s)g_1(r-s)\ dsdr\big )
=-f(t)+\int_0^t \int_0^r\pa k_2(0,s)f(r-s)\ dsdr.
$$
A standard integral equations argument now shows that $f\in L^2(0,T)$ implies $g_j\in L^2(0,T)$.

Now let
$\phi \in \mH^1_0$. Denoting the pairing $\mH^{-1}-\mH^1_0$ by $\< *,*\>$, we have 
\begin{eqnarray*}
\< u_t,\phi \> & = & \sum_j\int_0^{l_j}\phi_j(x)(u_j)_t(x,t)dx\nonumber\\
& = & \sum_j\int_0^{l_j}\phi_j(x)
\big ( g_j'(t-x)+\int_x^t  k_j(x,s)g_j'(t-s)ds\big )\nonumber\\
& = & \sum_j \int_0^{l_j}\phi_j'(x)
 g_j(t-x)dx +\phi (0)\sum_jg_j(t) \nonumber\\
&&+\sum_j \int_x^{l_j} \phi_j(x) \big ( k_j(x,x)g_j(t-x)+\int_x^t
(k_j)_s(x,s)g_j(t-s)ds\big )\ dx.\label{uno}
\end{eqnarray*}
This equation was derived assuming $f,g_j$ were regular, but 
by density it also holds for $f\in L^2(0,T).$
Inserting \rq{KN'} into\rq{uno}, we see that for $f\in L^2(0,T)$, the mapping
$t\mapsto \< u_t,\phi \>$ is continuous. Hence the mapping 
$t\mapsto  u_t$ is in $C(0,T; \mH^{-1}).$
This completes the proof of Theorem \ref{thm1} for small $t.$
For larger $t$, one must account for reflected waves. Although this causes more terms to appear in the formula for $\< u_t,\phi \>$, the formulas for reflected waves implicit in \rq{fold2} implies the continuity remains unchanged.}

		\section{Shape and velocity controllability on graphs}  \label{sv}
		
		In this section we prove the shape and velocity controllability on general graphs.  For readability, the proof and examples are given for the case of zero potentials. The corresponding equations for nonzero potentials are presented in the Appendix. 
		
		\subsection{Geometric preliminaries}
		
		In \cite{AZ}, shape/velocity controllability on a tree graph with Dirichlet boundary control was proved by representing a tree graph with a union of paths, and controlling the final shape/velocity on one path at a time. This approach was extended to general graphs in \cite{AZ2}, and we now adopt the approach and notation of that paper.

		Let $\vec{\Omega}=(V,\vec{E})$ be a directed graph obtained by orienting edges in $\Omega$.  We denote a directed edge $e_j$ from $v_i$ to $v_k$ as $e_j(v_i, v_k)$. When $e_j(v_i,v_k)$ is identified with the interval $(0,l_j)$, $v_i$ is identified with $x=0$ and $v_k$ is identified with $x=l_j$.  We denote the sets of indices for all outgoing edges of $v_i$ by $J^+(v_i)$ and for all incoming edges by $J^-(v_i)$. Since orienting $\Omega$ does not add or delete vertices or edges to/from the graph, we refer to the vertices and edges in $\vec{\Omega}$ by the same names as in $\Omega$.

		It was proven in \cite{AZ2} that for any graph $\Omega$,
		there is a directed acyclic graph (DAG), $\vec{\Omega}=(V,\vec{E})$ based on $\Omega$.
		For the rest of the text we assume all directed graphs are acyclic. See Figure \ref{D} for a DAG graph $\vec{\Omega}$ based on a graph $\Omega$.

		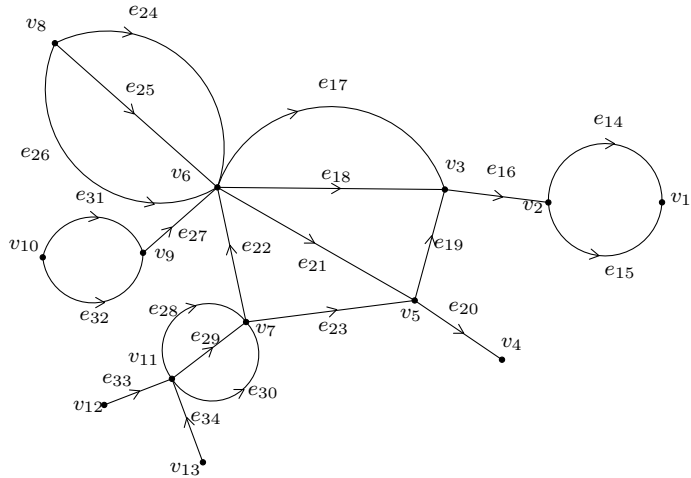
\begin{figure}[H]
			\begin{center}
				\begin{tikzpicture}[line cap=round,line join=round,>=triangle 45,x=1.0cm,y=1.0cm]
					\draw [line width=0.4pt] (-0.3825610335641724,6.0586212849287655)-- (1.778010245331328,4.140640290455626);
					\draw [line width=0.4pt] (0.7893020498686203,3.270395327795571)-- (1.778010245331328,4.140640290455626);
					\draw [line width=0.4pt] (1.778010245331328,4.140640290455626)-- (4.4004829655874715,2.635942828013578);
					\draw [line width=0.4pt] (2.159058366061286,2.349333787548426)-- (1.778010245331328,4.140640290455626);
					\draw [line width=0.4pt] (2.159058366061286,2.349333787548426)-- (4.4004829655874715,2.635942828013578);
					\draw [line width=0.4pt] (4.4004829655874715,2.635942828013578)-- (4.801735622238684,4.111979386409111);
					\draw [line width=0.4pt] (1.778010245331328,4.140640290455626)-- (4.801735622238684,4.111979386409111);
					\draw [shift={(3.285096116443921,3.6223556089478075)},line width=0.4pt]  plot[domain=0.31227216526678087:2.810363710063288,variable=\t]({1.*1.59371479066766*cos(\t r)+0.*1.59371479066766*sin(\t r)},{0.*1.59371479066766*cos(\t r)+1.*1.59371479066766*sin(\t r)});
					\draw [line width=0.4pt] (4.801735622238684,4.111979386409111)-- (6.177459016471415,3.94019817730482);
					\draw [line width=0.4pt] (4.4004829655874715,2.635942828013578)-- (5.561249579471338,1.8477679667344105);
					\draw [shift={(1.6642689638129242,1.9737947555669575)},line width=0.4pt]  plot[domain=0.6492284162694677:3.8006333935724523,variable=\t]({1.*0.6211651287048153*cos(\t r)+0.*0.6211651287048153*sin(\t r)},{0.*0.6211651287048153*cos(\t r)+1.*0.6211651287048153*sin(\t r)});
					\draw [shift={(1.70416128499514,1.921766880639986)},line width=0.4pt]  plot[domain=-2.587761257186781:0.7544377598491149,variable=\t]({1.*0.6242954542888733*cos(\t r)+0.*0.6242954542888733*sin(\t r)},{0.*0.6242954542888733*cos(\t r)+1.*0.6242954542888733*sin(\t r)});
					\draw [line width=0.4pt] (2.159058366061286,2.349333787548426)-- (1.1731882323724663,1.5934188463314032);
					\draw [line width=0.4pt] (1.1731882323724663,1.5934188463314032)-- (0.27331278288927857,1.2458889817575913);
					\draw [line width=0.4pt] (1.1731882323724663,1.5934188463314032)-- (1.5851489481736283,0.4827642591372834);
					\draw [shift={(0.3105618944053427,4.663498919290151)},line width=0.4pt]  plot[domain=-0.34228033026990623:2.0318956861262825,variable=\t]({1.*1.5578144332307726*cos(\t r)+0.*1.5578144332307726*sin(\t r)},{0.*1.5578144332307726*cos(\t r)+1.*1.5578144332307726*sin(\t r)});
					\draw [shift={(1.0063301312448432,5.447269391268226)},line width=0.4pt]  plot[domain=2.726941148346936:5.245859514689026,variable=\t]({1.*1.5174879919020363*cos(\t r)+0.*1.5174879919020363*sin(\t r)},{0.*1.5174879919020363*cos(\t r)+1.*1.5174879919020363*sin(\t r)});
					\draw [shift={(6.93308204983153,3.967974341690649)},line width=0.4pt]  plot[domain=-0.03519246039618107:3.1783353976437123,variable=\t]({1.*0.7561333770257241*cos(\t r)+0.*0.7561333770257241*sin(\t r)},{0.*0.7561333770257241*cos(\t r)+1.*0.7561333770257241*sin(\t r)});
					\draw [shift={(6.933068738038241,3.985147702909215)},line width=0.4pt]  plot[domain=3.201010385079611:6.2253178593475065,variable=\t]({1.*0.7569455140090021*cos(\t r)+0.*0.7569455140090021*sin(\t r)},{0.*0.7569455140090021*cos(\t r)+1.*0.7569455140090021*sin(\t r)});
					\draw [shift={(0.13101342621735887,3.053945165410281)},line width=0.4pt]  plot[domain=0.31767168727252215:2.9147675251604253,variable=\t]({1.*0.6929607383000072*cos(\t r)+0.*0.6929607383000072*sin(\t r)},{0.*0.6929607383000072*cos(\t r)+1.*0.6929607383000072*sin(\t r)});
					\draw [shift={(0.12111567660261073,3.271695656934737)},line width=0.4pt]  plot[domain=3.2343852675302767:6.281239252082257,variable=\t]({1.*0.6681876385224838*cos(\t r)+0.*0.6681876385224838*sin(\t r)},{0.*0.6681876385224838*cos(\t r)+1.*0.6681876385224838*sin(\t r)});
					\begin{scriptsize}
						\draw [fill=black] (7.688747236083874,3.9413696402539413) circle (1.0pt);
						\draw[color=black] (7.946249171179057,3.990787880592706) node {$v_1$};
						\draw [fill=black] (6.177459016471415,3.94019817730482) circle (1.0pt);
						\draw[color=black] (5.980572109488782,3.8823367323615185) node {$v_2$};
						\draw [fill=black] (4.801735622238684,4.111979386409111) circle (1.0pt);
						\draw[color=black] (4.950286201292499,4.424592473517456) node {$v_3$};
						\draw [fill=black] (5.561249579471338,1.8477679667344105) circle (1.0pt);
						\draw[color=black] (5.695887845381915,2.0522236059602306) node {$v_4$};
						\draw [fill=black] (4.4004829655874715,2.635942828013578) circle (1.0pt);
						\draw[color=black] (4.353804886020968,2.445359018298285) node {$v_5$};
						\draw [fill=black] (1.778010245331328,4.140640290455626) circle (1.0pt);
						\draw[color=black] (1.3036163420188156,4.30258493175737) node {$v_6$};
						\draw [fill=black] (2.159058366061286,2.349333787548426) circle (1.0pt);
						\draw[color=black] (2.4287970049173873,2.2420131153648084) node {$v_7$};
						\draw [fill=black] (-0.3825610335641724,6.0586212849287655) circle (1.0pt);
						\draw[color=black] (-0.6213915390847646,6.29537478050544) node {$v_8$};
						\draw [fill=black] (0.7893020498686203,3.270395327795571) circle (1.0pt);
						\draw[color=black] (1.0867140455564404,3.245186236503292) node {$v_9$};
						\draw [fill=black] (-0.5441973209342128,3.2097817200318057) circle (1.0pt);
						\draw[color=black] (-0.7840682614315461,3.367193778263378) node {$v_{10}$};
						\draw [fill=black] (1.1731882323724663,1.5934188463314032) circle (1.0pt);
						\draw[color=black] (0.8020297814495729,1.8624340965556523) node {$v_{11}$};
						\draw [fill=black] (0.27331278288927857,1.2458889817575913) circle (1.0pt);
						\draw[color=black] (0.06998453088905646,1.2252836006974261) node {$v_{12}$};
						\draw [fill=black] (1.5851489481736283,0.4827642591372834) circle (1.0pt);
						\draw[color=black] (1.3578419161344095,0.39834359543462194) node {$v_{13}$};
						\draw[color=black] (6.983745230627267,5.021073788788987) node {$e_{14}$};
						\draw[color=black] (7.1328655594451496,3.014727546512019) node {$e_{15}$};
						\draw[color=black] (5.56032391009293,4.356810505872964) node {$e_{16}$};
						\draw[color=black] (3.323518977824685,5.50910395582933) node {$e_{17}$};
						\draw[color=black] (3.3641881584113804,4.292374441161948) node {$e_{18}$};
						\draw[color=black] (4.860496691887921,3.3807501717922763) node {$e_{19}$};
						\draw[color=black] (5.058737349523687,2.5402537730005736) node {$e_{20}$};
						\draw[color=black] (3.052391107246716,3.109622301214308) node {$e_{21}$};
						\draw[color=black] (2.306789463157301,3.340080991205581) node {$e_{22}$};
						\draw[color=black] (3.323518977824685,2.2691259024226054) node {$e_{23}$};
						\draw[color=black] (0.8020297814495729,6.471607896381119) node {$e_{24}$}; 
						\draw[color=black] (0.7613606008628775,5.454878381713737) node {$e_{25}$};
						\draw[color=black] (-0.6349479326136631,4.587269195864237) node {$e_{26}$};
						\draw[color=black] (1.4578419161344095,3.5027577135523624) node {$e_{27}$};
						\draw[color=black] (1.0680524067298312,2.502930495347355) node {$e_{28}$};
						\draw[color=black] (1.6289697867123787,2.1420131153648084) node {$e_{29}$};
						\draw[color=black] (2.388127824330692,1.415073110102004) node {$e_{30}$};
						\draw[color=black] (0.09709731794685335,4.017900667650503) node {$e_{31}$};
						\draw[color=black] (0.15132289206244717,2.418246231240488) node {$e_{32}$};
						\draw[color=black] (0.4524507626404162,1.5726445871510744) node {$e_{33}$};
						\draw[color=black] (1.6289697867123787,1.0761632718795435) node {$e_{34}$};
						\draw [shift={(6.976967033862818,4.722833131153221)},rotate=30] (0,0) ++(0 pt,2.25pt) -- ++(1.9485571585149868pt,-3.375pt)--++(-3.8971143170299736pt,0 pt); 
						\draw [shift={(6.7736211309293415,3.2451862365032924)},rotate=20] (0,0) ++(0 pt,2.25pt) -- ++(1.9485571585149868pt,-3.375pt)--++(-3.8971143170299736pt,0 pt); 
						\draw [shift={(5.492791435233944,4.025689945077715)},rotate=20] (0,0) ++(0 pt,2.25pt) -- ++(1.9485571585149868pt,-3.375pt)--++(-3.8971143170299736pt,0 pt); 
						\draw [shift={(2.7668114349361033,5.1294414800604)},rotate=45] (0,0) ++(0 pt,2.25pt) -- ++(1.9485571585149868pt,-3.375pt)--++(-3.8971143170299736pt,0 pt); 
						\draw [shift={(3.334947435847778,4.1258825919673185)},rotate=30] (0,0) ++(0 pt,2.25pt) -- ++(1.9485571585149868pt,-3.375pt)--++(-3.8971143170299736pt,0 pt); 
						\draw [shift={(4.6134998028178975,3.419540479254073)},rotate=110] (0,0) ++(0 pt,2.25pt) -- ++(1.9485571585149868pt,-3.375pt)--++(-3.8971143170299736pt,0 pt) ; 
						\draw [shift={(4.9848210288134425,2.2391700690329808)},rotate=0] (0,0) ++(0 pt,2.25pt) -- ++(1.9485571585149868pt,-3.375pt)--++(-3.8971143170299736pt,0 pt); 
						\draw [shift={(3.001137206500698,3.438846132407628)},rotate=0] (0,0) ++(0 pt,2.25pt) -- ++(1.9485571585149868pt,-3.375pt)--++(-3.8971143170299736pt,0 pt); 
						\draw [shift={(3.2883930760596733,2.4937408479800154)},rotate=45] (0,0) ++(0 pt,2.25pt) -- ++(1.9485571585149868pt,-3.375pt)--++(-3.8971143170299736pt,0 pt) ; 
						\draw [shift={(0.5670521547348149,6.20005303637755)},rotate=20] (0,0) ++(0 pt,2.25pt) -- ++(1.9485571585149868pt,-3.375pt)--++(-3.8971143170299736pt,0 pt); 
						\draw [shift={(0.6140009260527342,5.17395392077809)},rotate=345] (0,0) ++(0 pt,2.25pt) -- ++(1.9485571585149868pt,-3.375pt)--++(-3.8971143170299736pt,0 pt) ; 
						\draw [shift={(0.863033552329616,3.936562306477112)},rotate=15] (0,0) ++(0 pt,2.25pt) -- ++(1.9485571585149868pt,-3.375pt)--++(-3.8971143170299736pt,0 pt); 
						\draw [shift={(1.1131075683935487,3.5554037080952514)},rotate=75] (0,0) ++(0 pt,2.25pt) -- ++(1.9485571585149868pt,-3.375pt)--++(-3.8971143170299736pt,0 pt); 
						\draw [shift={(1.958502367836447,3.2921471601192835)},rotate=130] (0,0) ++(0 pt,2.25pt) -- ++(1.9485571585149868pt,-3.375pt)--++(-3.8971143170299736pt,0 pt); 
						\draw [,shift={(1.4156426634275268,2.54303203463039)},rotate=45] (0,0) ++(0 pt,2.25pt) -- ++(1.9485571585149868pt,-3.375pt)--++(-3.8971143170299736pt,0 pt) ; 
						\draw [shift={(1.652948243036306,1.961274355760705)},rotate=70] (0,0) ++(0 pt,2.25pt) -- ++(1.9485571585149868pt,-3.375pt)--++(-3.8971143170299736pt,0 pt); 
						\draw [shift={(2.0641329657095455,1.411702959813001)},rotate=60] (0,0) ++(0 pt,2.25pt) -- ++(1.9485571585149868pt,-3.375pt)--++(-3.8971143170299736pt,0 pt);
						\draw [shift={(0.11743190824020103,3.746772797072534)},rotate=20] (0,0) ++(0 pt,2.25pt) -- ++(1.9485571585149868pt,-3.375pt)--++(-3.8971143170299736pt,0 pt); 
						\draw [shift={(0.19877026941359174,2.608035740645066)},rotate=30] (0,0) ++(0 pt,2.25pt) -- ++(1.9485571585149868pt,-3.375pt)--++(-3.8971143170299736pt,0 pt); 
						\draw [shift={(0.6725220297241353,1.400062690540079)},rotate=55] (0,0) ++(0 pt,2.25pt) -- ++(1.9485571585149868pt,-3.375pt)--++(-3.8971143170299736pt,0 pt); 
						\draw [shift={(1.3928157600996272,1.0012985013856708)},rotate=135] (0,0) ++(0 pt,2.25pt) -- ++(1.9485571585149868pt,-3.375pt)--++(-3.8971143170299736pt,0 pt); 
					\end{scriptsize}
				\end{tikzpicture}
				\caption{A DAG $\vec{\Omega}$ based on $\Omega$. } \label{D}
			\end{center}
		\end{figure}
		
		In a directed graph $\vec{\Omega}$, a \textbf{source} is a vertex without incoming edges; a \textbf{sink} is a vertex without outgoing edges. We denote the sets of sources and sinks by $\vec{\Omega}^+$ and $\vec{\Omega}^-$ respectively.  
		
		The following lemma, proven in \cite{AZ2}, ensures a linear {\bf what does "linear" order mean?} ordering of vertices that is consistent with edge directions exists.  Such a linear ordering is critical in solving the control problem.

		\begin{lemma} \label{topoorder}
			Any DAG $\vec{\Omega}$ has at least one linear ordering $W(V, <)$ of the set $V$ such that for every directed path from $v_i$ to $v_k$, $v_k < v_i$ {\bf should this be $k<i$?}in the ordering. 
		\end{lemma}

		\begin{definition} \label{st}
			
			Let $\vec{\Omega}$ be a DAG of $\Omega$.  Let $U$ be a union of directed paths.  We say $U$ is a  \textbf{tangle-free (TF) path union} of $\Omega$ if $U$ satisfies the conditions: 
			
			\begin{enumerate} 
				
				\item The direction of all edges are the same as the direction of all paths they are on. 
				
				\item All paths are disjoint except for the starting and finishing vertices.						
				
				\item If a finishing vertex $v$ of a path is the starting vertex of another path, there must be an incoming edge of $v$ that is not a finishing edge, and an outgoing edge of $v$ that is not a starting edge. 
				
				\item $\vec{\Omega}=\cup_{P \in U} P$.
				
			\end{enumerate}
		\end{definition}

		\begin{definition}
			We say is $I^*$ is the \textbf{single-track (ST) active set of vertex indices} of $\vec{\Omega}$ if
			\begin{equation} \label{vstar}I^{*}=\{i: v_i \in \Omega^+\}, \end{equation} 
			$J^*$ is a \textbf{ST active set of edge indices} if for all \je{$v_i \in I\setminus \Gamma$,
			\begin{equation} \label{jstar}J^*(v_i)= \textrm{all but one elements of } J^+(v_i) , \end{equation} 
			and for $v_i\in \Gamma ,$
			\begin{equation} \label{jstarb}\je{J^*(v_i)=  J^+(v_i) . }\end{equation} }
			The set $\{I^*, J^*\}$ is then called a \textbf{ST active set}.  
		\end{definition}
		\begin{example} \label{egst}
		\je{	For the DAG $\vec{\Omega}$ in Figure \ref{D}, a ST active set is $V^*=\{v_8, v_{10}\}$, $J^*(v_1)=\emptyset$,  $J^*(v_2)=\{15\}$,  $J^*(v_3)=\emptyset$,  $J^*(v_4)=\emptyset$,  $J^*(v_5)=\{19\}$,  $J^*(v_6)=\{18,21\}$,  $J^*(v_7)=\{23\}$,  $J^*(v_8)=\{25,26\}$,  $J^*(v_9)=\emptyset$,  $J^*(v_{10})=\{32\}$,  $J^*(v_{11})=\{29,30\}$,  $J^*(v_{12})=\{33\}$,  $J^*(v_{13})=\{ 34\}$.}  See Figure \ref{activeset}. 
			
			One ST path union associated with $\{I^*, J^*\}$ is then $P_1=(v_8, e_{24}, v_6, e_{17}, v_3, e_{16},v_2, e_{14}, v_1)$, $P_2=(v_8, e_{26}, v_6)$, $P_3=(v_8, e_{25}, v_6)$, $P_4=(v_6, e_{18}, v_3)$, $P_5=(v_5, e_{19}, v_3)$, $P_6=(v_2, e_{15}, v_1)$, $P_7=( v_{10}, e_{31}, v_9, e_{27}, v_6)$, $P_8=(v_{10}, e_{32}, v_9)$, $P_9=(v_{12}, e_{33}, v_{11},e_{28}, v_7, e_{22}, v_6)$, $P_{10}=(v_{13}, e_{34}, v_{11})$, $P_{11}=(v_{11}, e_{29}, v_7)$, $P_{12}=(v_{11}, e_{30}, v_7)$, $P_{13}=(v_6, e_{21}, v_5,e_{20}, v_4)$, $P_{14}=(v_7, e_{23}, v_5)$.
			
			\je{Note that $|U|=|I^*|+|J^*|=14$. }
		\end{example}
		
		
		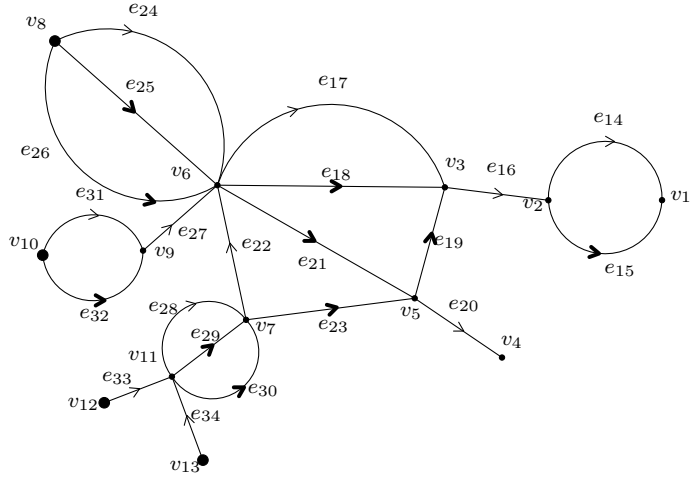
\begin{figure}[H]
			\begin{center}
				\begin{tikzpicture}[line cap=round,line join=round,>=triangle 45,x=1.0cm,y=1.0cm]
					\draw [line width=0.4pt] (-0.3825610335641724,6.0586212849287655)-- (1.778010245331328,4.140640290455626);
					\draw [line width=0.4pt] (0.7893020498686203,3.270395327795571)-- (1.778010245331328,4.140640290455626);
					\draw [line width=0.4pt] (1.778010245331328,4.140640290455626)-- (4.4004829655874715,2.635942828013578);
					\draw [line width=0.4pt] (2.159058366061286,2.349333787548426)-- (1.778010245331328,4.140640290455626);
					\draw [line width=0.4pt] (2.159058366061286,2.349333787548426)-- (4.4004829655874715,2.635942828013578);
					\draw [line width=0.4pt] (4.4004829655874715,2.635942828013578)-- (4.801735622238684,4.111979386409111);
					\draw [line width=0.4pt] (1.778010245331328,4.140640290455626)-- (4.801735622238684,4.111979386409111);
					\draw [shift={(3.285096116443921,3.6223556089478075)},line width=0.4pt]  plot[domain=0.31227216526678087:2.810363710063288,variable=\t]({1.*1.59371479066766*cos(\t r)+0.*1.59371479066766*sin(\t r)},{0.*1.59371479066766*cos(\t r)+1.*1.59371479066766*sin(\t r)});
					\draw [line width=0.4pt] (4.801735622238684,4.111979386409111)-- (6.177459016471415,3.94019817730482);
					\draw [line width=0.4pt] (4.4004829655874715,2.635942828013578)-- (5.561249579471338,1.8477679667344105);
					\draw [shift={(1.6642689638129242,1.9737947555669575)},line width=0.4pt]  plot[domain=0.6492284162694677:3.8006333935724523,variable=\t]({1.*0.6211651287048153*cos(\t r)+0.*0.6211651287048153*sin(\t r)},{0.*0.6211651287048153*cos(\t r)+1.*0.6211651287048153*sin(\t r)});
					\draw [shift={(1.70416128499514,1.921766880639986)},line width=0.4pt]  plot[domain=-2.587761257186781:0.7544377598491149,variable=\t]({1.*0.6242954542888733*cos(\t r)+0.*0.6242954542888733*sin(\t r)},{0.*0.6242954542888733*cos(\t r)+1.*0.6242954542888733*sin(\t r)});
					\draw [line width=0.4pt] (2.159058366061286,2.349333787548426)-- (1.1731882323724663,1.5934188463314032);
					\draw [line width=0.4pt] (1.1731882323724663,1.5934188463314032)-- (0.27331278288927857,1.2458889817575913);
					\draw [line width=0.4pt] (1.1731882323724663,1.5934188463314032)-- (1.5851489481736283,0.4827642591372834);
					\draw [shift={(0.3105618944053427,4.663498919290151)},line width=0.4pt]  plot[domain=-0.34228033026990623:2.0318956861262825,variable=\t]({1.*1.5578144332307726*cos(\t r)+0.*1.5578144332307726*sin(\t r)},{0.*1.5578144332307726*cos(\t r)+1.*1.5578144332307726*sin(\t r)});
					\draw [shift={(1.0063301312448432,5.447269391268226)},line width=0.4pt]  plot[domain=2.726941148346936:5.245859514689026,variable=\t]({1.*1.5174879919020363*cos(\t r)+0.*1.5174879919020363*sin(\t r)},{0.*1.5174879919020363*cos(\t r)+1.*1.5174879919020363*sin(\t r)});
					\draw [shift={(6.93308204983153,3.967974341690649)},line width=0.4pt]  plot[domain=-0.03519246039618107:3.1783353976437123,variable=\t]({1.*0.7561333770257241*cos(\t r)+0.*0.7561333770257241*sin(\t r)},{0.*0.7561333770257241*cos(\t r)+1.*0.7561333770257241*sin(\t r)});
					\draw [shift={(6.933068738038241,3.985147702909215)},line width=0.4pt]  plot[domain=3.201010385079611:6.2253178593475065,variable=\t]({1.*0.7569455140090021*cos(\t r)+0.*0.7569455140090021*sin(\t r)},{0.*0.7569455140090021*cos(\t r)+1.*0.7569455140090021*sin(\t r)});
					\draw [shift={(0.13101342621735887,3.053945165410281)},line width=0.4pt]  plot[domain=0.31767168727252215:2.9147675251604253,variable=\t]({1.*0.6929607383000072*cos(\t r)+0.*0.6929607383000072*sin(\t r)},{0.*0.6929607383000072*cos(\t r)+1.*0.6929607383000072*sin(\t r)});
					\draw [shift={(0.12111567660261073,3.271695656934737)},line width=0.4pt]  plot[domain=3.2343852675302767:6.281239252082257,variable=\t]({1.*0.6681876385224838*cos(\t r)+0.*0.6681876385224838*sin(\t r)},{0.*0.6681876385224838*cos(\t r)+1.*0.6681876385224838*sin(\t r)});
					\begin{scriptsize}
						\draw [fill=black] (7.688747236083874,3.9413696402539413) circle (1.0pt);
						\draw[color=black] (7.946249171179057,3.990787880592706) node {$v_1$};
						\draw [fill=black] (6.177459016471415,3.94019817730482) circle (1.0pt);
						\draw[color=black] (5.980572109488782,3.8823367323615185) node {$v_2$};
						\draw [fill=black] (4.801735622238684,4.111979386409111) circle (1.0pt);
						\draw[color=black] (4.950286201292499,4.424592473517456) node {$v_3$};
						\draw [fill=black] (5.561249579471338,1.8477679667344105) circle (1.0pt);
						\draw[color=black] (5.695887845381915,2.0522236059602306) node {$v_4$};
						\draw [fill=black] (4.4004829655874715,2.635942828013578) circle (1.0pt);
						\draw[color=black] (4.353804886020968,2.445359018298285) node {$v_5$};
						\draw [fill=black] (1.778010245331328,4.140640290455626) circle (1.0pt);
						\draw[color=black] (1.3036163420188156,4.30258493175737) node {$v_6$};
						\draw [fill=black] (2.159058366061286,2.349333787548426) circle (1.0pt);
						\draw[color=black] (2.4287970049173873,2.2420131153648084) node {$v_7$};
						\draw [fill=black] (-0.3825610335641724,6.0586212849287655) circle (2.0pt);
						\draw[color=black] (-0.6213915390847646,6.29537478050544) node {$v_8$};
						\draw [fill=black] (0.7893020498686203,3.270395327795571) circle (1.0pt);
						\draw[color=black] (1.0867140455564404,3.245186236503292) node {$v_9$};
						\draw [fill=black] (-0.5441973209342128,3.2097817200318057) circle (2.0pt);
						\draw[color=black] (-0.7840682614315461,3.367193778263378) node {$v_{10}$};
						\draw [fill=black] (1.1731882323724663,1.5934188463314032) circle (1.0pt);
						\draw[color=black] (0.8020297814495729,1.8624340965556523) node {$v_{11}$};
						\draw [fill=black] (0.27331278288927857,1.2458889817575913) circle (2.0pt);
						\draw[color=black] (0,1.2252836006974261) node {$v_{12}$};
						\draw [fill=black] (1.5851489481736283,0.4827642591372834) circle (2.0pt);
						\draw[color=black] (1.3578419161344095,0.39834359543462194) node {$v_{13}$};
						\draw[color=black] (6.983745230627267,5.021073788788987) node {$e_{14}$};
						\draw[color=black] (7.1328655594451496,3.014727546512019) node {$e_{15}$};
						\draw[color=black] (5.56032391009293,4.356810505872964) node {$e_{16}$};
						\draw[color=black] (3.323518977824685,5.50910395582933) node {$e_{17}$};
						\draw[color=black] (3.3641881584113804,4.292374441161948) node {$e_{18}$};
						\draw[color=black] (4.860496691887921,3.3807501717922763) node {$e_{19}$};
						\draw[color=black] (5.058737349523687,2.5402537730005736) node {$e_{20}$};
						\draw[color=black] (3.052391107246716,3.109622301214308) node {$e_{21}$};
						\draw[color=black] (2.306789463157301,3.340080991205581) node {$e_{22}$};
						\draw[color=black] (3.323518977824685,2.2691259024226054) node {$e_{23}$};
						\draw[color=black] (0.8020297814495729,6.471607896381119) node {$e_{24}$}; 
						\draw[color=black] (0.7613606008628775,5.454878381713737) node {$e_{25}$};
						\draw[color=black] (-0.6349479326136631,4.587269195864237) node {$e_{26}$};
						\draw[color=black] (1.4578419161344095,3.5027577135523624) node {$e_{27}$};
						\draw[color=black] (1.0680524067298312,2.502930495347355) node {$e_{28}$};
						\draw[color=black] (1.6289697867123787,2.1420131153648084) node {$e_{29}$};
						\draw[color=black] (2.388127824330692,1.415073110102004) node {$e_{30}$};
						\draw[color=black] (0.09709731794685335,4.017900667650503) node {$e_{31}$};
						\draw[color=black] (0.15132289206244717,2.418246231240488) node {$e_{32}$};
						\draw[color=black] (0.4524507626404162,1.5726445871510744) node {$e_{33}$};
						\draw[color=black] (1.6289697867123787,1.0761632718795435) node {$e_{34}$};
						\draw [shift={(6.976967033862818,4.722833131153221)},rotate=30] (0,0) ++(0 pt,2.25pt) -- ++(1.9485571585149868pt,-3.375pt)--++(-3.8971143170299736pt,0 pt); 
						\draw [line width=1.5pt,shift={(6.7736211309293415,3.2451862365032924)},rotate=20] (0,0) ++(0 pt,2.25pt) -- ++(1.9485571585149868pt,-3.375pt)--++(-3.8971143170299736pt,0 pt); 
						\draw [shift={(5.492791435233944,4.025689945077715)},rotate=20] (0,0) ++(0 pt,2.25pt) -- ++(1.9485571585149868pt,-3.375pt)--++(-3.8971143170299736pt,0 pt); 
						\draw [shift={(2.7668114349361033,5.1294414800604)},rotate=45] (0,0) ++(0 pt,2.25pt) -- ++(1.9485571585149868pt,-3.375pt)--++(-3.8971143170299736pt,0 pt); 
						\draw [line width=1.5pt,shift={(3.334947435847778,4.1258825919673185)},rotate=30] (0,0) ++(0 pt,2.25pt) -- ++(1.9485571585149868pt,-3.375pt)--++(-3.8971143170299736pt,0 pt); 
						\draw [line width=1.5pt,shift={(4.6134998028178975,3.419540479254073)},rotate=110] (0,0) ++(0 pt,2.25pt) -- ++(1.9485571585149868pt,-3.375pt)--++(-3.8971143170299736pt,0 pt) ; 
						\draw [shift={(4.9848210288134425,2.2391700690329808)},rotate=0] (0,0) ++(0 pt,2.25pt) -- ++(1.9485571585149868pt,-3.375pt)--++(-3.8971143170299736pt,0 pt); 
						\draw [line width=1.5pt,shift={(3.001137206500698,3.438846132407628)},rotate=0] (0,0) ++(0 pt,2.25pt) -- ++(1.9485571585149868pt,-3.375pt)--++(-3.8971143170299736pt,0 pt); 
						\draw [line width=1.5pt,shift={(3.2883930760596733,2.4937408479800154)},rotate=45] (0,0) ++(0 pt,2.25pt) -- ++(1.9485571585149868pt,-3.375pt)--++(-3.8971143170299736pt,0 pt) ; 
						\draw [shift={(0.5670521547348149,6.20005303637755)},rotate=20] (0,0) ++(0 pt,2.25pt) -- ++(1.9485571585149868pt,-3.375pt)--++(-3.8971143170299736pt,0 pt); 
						\draw [line width=1.5pt,shift={(0.6140009260527342,5.17395392077809)},rotate=345] (0,0) ++(0 pt,2.25pt) -- ++(1.9485571585149868pt,-3.375pt)--++(-3.8971143170299736pt,0 pt) ; 
						\draw [line width=1.5pt,shift={(0.863033552329616,3.936562306477112)},rotate=15] (0,0) ++(0 pt,2.25pt) -- ++(1.9485571585149868pt,-3.375pt)--++(-3.8971143170299736pt,0 pt); 
						\draw [shift={(1.1131075683935487,3.5554037080952514)},rotate=75] (0,0) ++(0 pt,2.25pt) -- ++(1.9485571585149868pt,-3.375pt)--++(-3.8971143170299736pt,0 pt); 
						\draw [shift={(1.958502367836447,3.2921471601192835)},rotate=130] (0,0) ++(0 pt,2.25pt) -- ++(1.9485571585149868pt,-3.375pt)--++(-3.8971143170299736pt,0 pt); 
						\draw [,shift={(1.4156426634275268,2.54303203463039)},rotate=45] (0,0) ++(0 pt,2.25pt) -- ++(1.9485571585149868pt,-3.375pt)--++(-3.8971143170299736pt,0 pt) ; 
						\draw [line width=1.5pt,shift={(1.652948243036306,1.961274355760705)},rotate=70] (0,0) ++(0 pt,2.25pt) -- ++(1.9485571585149868pt,-3.375pt)--++(-3.8971143170299736pt,0 pt); 
						\draw [line width=1.5pt,shift={(2.0641329657095455,1.411702959813001)},rotate=60] (0,0) ++(0 pt,2.25pt) -- ++(1.9485571585149868pt,-3.375pt)--++(-3.8971143170299736pt,0 pt);
						\draw [shift={(0.11743190824020103,3.746772797072534)},rotate=20] (0,0) ++(0 pt,2.25pt) -- ++(1.9485571585149868pt,-3.375pt)--++(-3.8971143170299736pt,0 pt); 
						\draw [line width=1.5pt,shift={(0.19877026941359174,2.608035740645066)},rotate=30] (0,0) ++(0 pt,2.25pt) -- ++(1.9485571585149868pt,-3.375pt)--++(-3.8971143170299736pt,0 pt); 
						\draw [shift={(0.6725220297241353,1.400062690540079)},rotate=55] (0,0) ++(0 pt,2.25pt) -- ++(1.9485571585149868pt,-3.375pt)--++(-3.8971143170299736pt,0 pt); 
						\draw [shift={(1.3928157600996272,1.0012985013856708)},rotate=135] (0,0) ++(0 pt,2.25pt) -- ++(1.9485571585149868pt,-3.375pt)--++(-3.8971143170299736pt,0 pt); 
					\end{scriptsize}
				\end{tikzpicture}
				\caption{$\vec{\Omega}$ with ST active vertices and edges emphasized.} \label{activeset}
			\end{center}
		\end{figure}

		We list some facts  about a DAG of $\Omega$ for which $\{ I^*,J^*\}$ is a ST active set.
		\je{These facts were proven for a slightly different setting in  \cite{AZ2}, but are easily seen to still apply here. }

1. For each $i\in I^*$, there will be $|J(v_i)|-1$ active edges incident to $v_i$.

2. There is a bijection between controls and heads of paths in $U$. In particular, at every head of path, there will either be an edge control or a vertex control. 
		
		\begin{definition}
			Let $U$ be a TF path union of a DAG $\vec{\Omega}$.  We define \textbf{depth} of an edge $e_j(v_i, v_k)$ in $\vec{\Omega}, U$, recursively as 
			
			\begin{equation} \label{depthsdef}
				\operatorname{depth}(e_j(v_i, v_k))=\begin{cases}l_j, & v_k {\rm \ is \ a \ finishing \ edge \ in \ }  U,\\ l_j + \max_{r \in J^+(v_k)} \operatorname{depth}(e_r), & {\rm otherwise}.\end{cases} 
			\end{equation}
		\end{definition}
		Depth of an edge depends on the choice of orientation and TF path union.  
		
		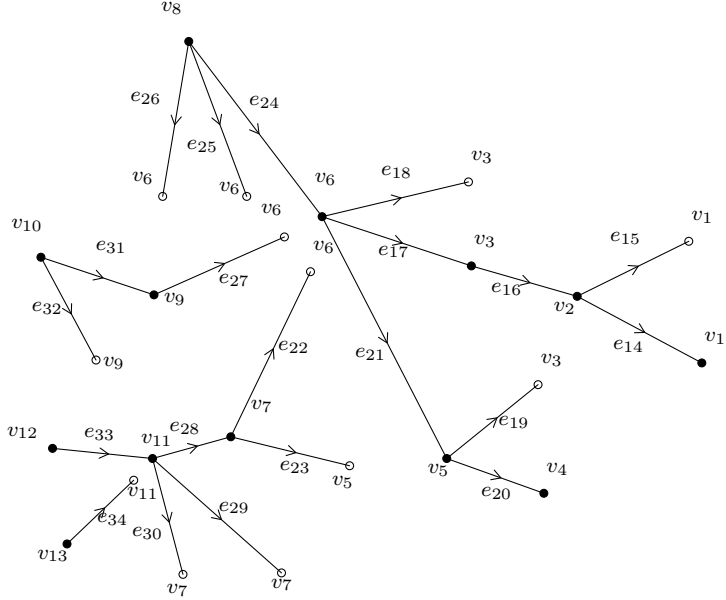
\begin{figure}[H]
			\begin{center}
				\begin{tikzpicture}[line cap=round,line join=round,>=triangle 45,x=1.0cm,y=1.0cm] 
					\draw [line width=0.4pt] (4.557245959441893,4.466660207436372)--(6.214226985575773,3.580368030667085);
					\draw [line width=0.4pt] (2.823196048371554,2.3087314292155097)-- (4.114099871057251,1.8463181195967533);
					\draw [line width=0.4pt] (-0.6063693313008958,7.857691144640581)-- (1.1662150222376737,5.526357375312685);
					\draw [line width=0.4pt] (-2.5716258971806125,4.986875180757469)-- (-1.068782640919653,4.485927428670482);
					\draw [line width=0.4pt] (-1.088049862153767,2.3087314292155052)-- (-0.0476199155115637,2.597739747727228);
					\draw [line width=0.4pt] (-2.417488127307694,2.443601977854309)-- (-1.088049862153767,2.3087314292155052);
					\draw [line width=0.4pt] (1.1662150222376737,5.526357375312685)-- (3.150738809351506,4.87127185335278);
					\draw [line width=0.4pt] (3.150738809351506,4.87127185335278)-- (4.557245959441893,4.466660207436372);
					\draw [line width=0.4pt] (-0.6063693313008958,7.857691144640581)-- (-0.9531793135149637,5.796098472590292);
					\draw [line width=0.4pt] (-0.6063693313008958,7.857691144640581)-- (0.16431951806369943,5.796098472590292);
					\draw [line width=0.4pt] (1.1662150222376737,5.526357375312685)-- (3.1122043668832764,5.988770684931441);
					\draw [line width=0.4pt] (4.557245959441893,4.466660207436372)-- (6.040821994468739,5.198814614332732);
					\draw [line width=0.4pt] (2.823196048371554,2.3087314292155097)-- (4.037030986120791,3.291359712155362);
					\draw [line width=0.4pt] (-0.0476199155115637,2.597739747727228)-- (1.532292225685856,2.212395323044931);
					\draw [line width=0.4pt] (-2.5716258971806125,4.986875180757469)-- (-1.8394714902842482,3.6189024731353148);
					\draw [line width=0.4pt] (-2.2248159149665447,1.1719653764027291)-- (-1.3385237381972612,2.0197231107037825);
					\draw [line width=0.4pt] (-0.0476199155115637,2.597739747727228)-- (1.0120772523647543,4.794202968416321);
					\draw [line width=0.4pt] (-1.068782640919653,4.485927428670482)-- (0.6652672701506863,5.256616278035077);
					\draw [line width=0.4pt] (1.1662150222376737,5.526357375312685)-- (2.823196048371554,2.3087314292155097);
					\draw [line width=0.4pt] (-1.088049862153767,2.3087314292155052)-- (-0.6834382162373547,0.7673537304863175);
					\draw [line width=0.4pt] (-1.088049862153767,2.3087314292155052)-- (0.626732827682457,0.7866209517204323);
					\draw [shift={(	5.385736473,4.023514119	)},rotate=	1.858398768	] (0,0) ++(0 pt,2.25pt) -- ++(1.9485571585149868pt,-3.375pt)--++(-3.8971143170299736pt,0 pt);
					\draw [shift={(	3.46864796,2.077524774	)},rotate=	10.29202116	] (0,0) ++(0 pt,2.25pt) -- ++(1.9485571585149868pt,-3.375pt)--++(-3.8971143170299736pt,0 pt);
					\draw [shift={(	0.279922845,6.69202426	)},rotate=	-22.75316187	] (0,0) ++(0 pt,2.25pt) -- ++(1.9485571585149868pt,-3.375pt)--++(-3.8971143170299736pt,0 pt);
					\draw [shift={(	-1.820204269,4.736401305	)},rotate=	11.56505118	] (0,0) ++(0 pt,2.25pt) -- ++(1.9485571585149868pt,-3.375pt)--++(-3.8971143170299736pt,0 pt);
					\draw [shift={(	-0.567834889,2.453235588	)},rotate=	45.524111	] (0,0) ++(0 pt,2.25pt) -- ++(1.9485571585149868pt,-3.375pt)--++(-3.8971143170299736pt,0 pt);
					\draw [shift={(	-1.752768995,2.376166704	)},rotate=	24.2072035	] (0,0) ++(0 pt,2.25pt) -- ++(1.9485571585149868pt,-3.375pt)--++(-3.8971143170299736pt,0 pt);
					\draw [shift={(	2.158476916,5.198814614	)},rotate=	11.73209379	] (0,0) ++(0 pt,2.25pt) -- ++(1.9485571585149868pt,-3.375pt)--++(-3.8971143170299736pt,0 pt);
					\draw [shift={(	3.853992384,4.66896603	)},rotate=	13.95099521	] (0,0) ++(0 pt,2.25pt) -- ++(1.9485571585149868pt,-3.375pt)--++(-3.8971143170299736pt,0 pt);
					\draw [shift={(	-0.779774322,6.826894809	)},rotate=	-70	] (0,0) ++(0 pt,2.25pt) -- ++(1.9485571585149868pt,-3.375pt)--++(-3.8971143170299736pt,0 pt);
					\draw [shift={(	-0.221024907,6.826894809	)},rotate=	-39.50265945	] (0,0) ++(0 pt,2.25pt) -- ++(1.9485571585149868pt,-3.375pt)--++(-3.8971143170299736pt,0 pt);
					\draw [shift={(	2.139209695,5.75756403	)},rotate=	43.3669307	] (0,0) ++(0 pt,2.25pt) -- ++(1.9485571585149868pt,-3.375pt)--++(-3.8971143170299736pt,0 pt);
					\draw [shift={(	5.299033977,4.832737411	)},rotate=	56.26663836	] (0,0) ++(0 pt,2.25pt) -- ++(1.9485571585149868pt,-3.375pt)--++(-3.8971143170299736pt,0 pt);
					\draw [shift={(	3.430113517,2.800045571	)},rotate=	68.99099404	] (0,0) ++(0 pt,2.25pt) -- ++(1.9485571585149868pt,-3.375pt)--++(-3.8971143170299736pt,0 pt);
					\draw [shift={(	0.742336155,2.405067535	)},rotate=	16.293039	] (0,0) ++(0 pt,2.25pt) -- ++(1.9485571585149868pt,-3.375pt)--++(-3.8971143170299736pt,0 pt);
					\draw [shift={(	-2.205548694,4.302888827	)},rotate=	-31.84380855	] (0,0) ++(0 pt,2.25pt) -- ++(1.9485571585149868pt,-3.375pt)--++(-3.8971143170299736pt,0 pt);
					\draw [shift={(	-1.781669827,1.595844244	)},rotate=	73.72696998	] (0,0) ++(0 pt,2.25pt) -- ++(1.9485571585149868pt,-3.375pt)--++(-3.8971143170299736pt,0 pt);
					\draw [shift={(	-0.885744039,	1.53804258	)},rotate=	-45.2916961	] (0,0) ++(0 pt,2.25pt) -- ++(1.9485571585149868pt,-3.375pt)--++(-3.8971143170299736pt,0 pt);
					\draw [shift={(	-0.230658517,	1.54767619	)},rotate=	-11.59355625	] (0,0) ++(0 pt,2.25pt) -- ++(1.9485571585149868pt,-3.375pt)--++(-3.8971143170299736pt,0 pt);
					\draw [shift={(	0.482228668,3.695971358	)},rotate=	94.24472924	] (0,0) ++(0 pt,2.25pt) -- ++(1.9485571585149868pt,-3.375pt)--++(-3.8971143170299736pt,0 pt);
					\draw [shift={(	-0.201757685,4.871271853	)},rotate=	53.96248897	] (0,0) ++(0 pt,2.25pt) -- ++(1.9485571585149868pt,-3.375pt)--++(-3.8971143170299736pt,0 pt);
					\draw [shift={(	1.994705535,3.917544402	)},rotate=	-32.75288843	] (0,0) ++(0 pt,2.25pt) -- ++(1.9485571585149868pt,-3.375pt)--++(-3.8971143170299736pt,0 pt);
					\begin{scriptsize}
						\draw [fill=black] (4.557245959441893,4.466660207436372) circle (1.5pt);
						\draw[color=black] (4.393474578951916,4.312522437563449) node {$v_2$};
						\draw [fill=black] (6.214226985575773,3.580368030667085) circle (1.5pt);
						\draw[color=black] (6.397265587299864,3.9079107916470375) node {$v_1$};
						\draw[color=black] (5.241232313252971,3.8115746854764634) node {$e_{14}$};
						\draw [fill=black] (3.150738809351506,4.87127185335278) circle (1.5pt);
						\draw[color=black] (3.333777411075598,5.198814614332732) node {$v_3$};
						\draw [fill=black] (1.1662150222376737,5.526357375312685) circle (1.5pt);
						\draw[color=black] (1.214383075322961,6.0273051273996705) node {$v_6$};
						\draw [fill=black] (2.823196048371554,2.3087314292155097) circle (1.5pt);
						\draw[color=black] (2.717226331583922,2.1545936593425865) node {$v_5$};
						\draw [fill=black] (4.114099871057251,1.8463181195967533) circle (1.5pt);
						\draw[color=black] (4.297138472781342,2.1738608805767012) node {$v_4$};
						\draw[color=black] (3.487915180948517,1.8463181195967489) node {$e_{20}$};
						\draw [fill=black] (-0.6063693313008958,7.857691144640581) circle (1.5pt);
						\draw[color=black] (-0.8279423754932165,8.300837233025224) node {$v_8$};
						\draw[color=black] (0.40515978349013604,7.048467852807757) node {$e_{24}$};
						\draw [fill=black] (-1.068782640919653,4.485927428670482) circle (1.5pt);
						\draw[color=black] (-0.7894079330249866,4.428125764968138) node {$v_9$};
						\draw [fill=black] (-2.5716258971806125,4.986875180757469) circle (1.5pt);
						\draw[color=black] (-2.754664498904704,5.410754047907995) node {$v_{10}$};
						\draw[color=black] (-1.6371656673260415,5.121745729396272) node {$e_{31}$};
						\draw [fill=black] (-0.0476199155115637,2.597739747727228) circle (1.5pt);
						\draw[color=black] (0.36662534102190625,3.00235139364364) node {$v_7$};
						\draw [fill=black] (-1.088049862153767,2.3087314292155052) circle (1.5pt);
						\draw[color=black] (-1.0279423754932165,2.558257553172012) node {$v_{11}$};
						\draw[color=black] (-0.6545373843861825,2.713343075131917) node {$e_{28}$};
						\draw [] (-0.6834382162373547,0.7673537304863175) circle (1.5pt);
						\draw[] (-0.7508734905567569,0.5554142969110542) node {$v_{7}$};
						\draw [] (0.626732827682457,0.7866209517204323) circle (1.5pt);
						\draw[color=black] (0.6363664382995146,0.6324831818475135) node {$v_{7}$};
						\draw [fill=black] (-2.417488127307694,2.443601977854309) circle (1.5pt);
						\draw[color=black] (-2.812466162607049,2.6555414114295726) node {$v_{12}$};
						\draw[color=black] (-1.7913034371989605,2.6362741901954574) node {$e_{33}$};
						\draw [fill=black] (-2.2248159149665447,1.1719653764027291) circle (1.5pt);
						\draw[color=black] (-2.4078545166906364,1.0178276065298104) node {$v_{13}$};
						\draw[color=black] (2.1199424733263603,5.063944065693928) node {$e_{17}$};
						\draw[color=black] (3.622785729587321,4.562996313606942) node {$e_{16}$};
						\draw [] (-0.9531793135149637,5.796098472590292) circle (1.5pt);
						\draw[color=black] (-1.213286800175514,6.046572348633785) node {$v_6$};
						\draw[color=black] (-1.1747523577072843,7.1062695165101015) node {$e_{26}$};
						\draw [] (0.16431951806369943,5.796098472590292) circle (1.5pt);
						\draw[color=black] (-0.05725352612862117,5.911701799994981) node {$v_6$};
						\draw[color=black] (-0.42333072957680395,6.470451215784312) node {$e_{25}$};
						\draw [] (3.1122043668832764,5.988770684931441) circle (1.5pt);
						\draw[color=black] (3.295242968607368,6.316313445911393) node {$v_3$};
						\draw[color=black] (2.1584769157945902,6.123641233570244) node {$e_{18}$};
						\draw [] (6.040821994468739,5.198814614332732) circle (1.5pt);
						\draw[color=black] (6.22386059619283,5.526357375312684) node {$v_1$};
						\draw[color=black] (5.202697870784742,5.237349056800961) node {$e_{15}$};
						\draw [] (4.037030986120791,3.291359712155362) circle (1.5pt);
						\draw[color=black] (4.220069587844883,3.6189024731353148) node {$v_3$};
						\draw[color=black] (3.7191218357578957,2.7904119600683766) node {$e_{19}$};
						\draw [] (1.532292225685856,2.212395323044931) circle (1.5pt);
						\draw[color=black] (1.4455897301323397,2.0004558894696673) node {$v_5$};
						\draw[color=black] (0.8097714294065486,2.193128101810816) node {$e_{23}$};
						\draw [] (-1.8394714902842482,3.6189024731353148) circle (1.5pt);
						\draw[color=black] (-1.5986312248578116,3.5803680306670853) node {$v_9$};
						\draw[color=black] (-2.484923401627096,4.312522437563449) node {$e_{32}$};
						\draw [] (-1.3385237381972612,2.0197231107037825) circle (1.5pt);
						\draw[color=black] (-1.2132868001755142,1.8463181195967489) node {$v_{11}$};
						\draw[color=black] (-1.6178984460919266,1.4995081373826817) node {$e_{34}$};
						\draw[color=black] (-1.15358963229919556,1.3187753586167965) node {$e_{30}$};
						\draw[color=black] (0,1.6573098010850261) node {$e_{29}$};
						\draw [] (1.0120772523647543,4.794202968416321) circle (1.5pt);
						\draw[color=black] (1.1951158540888462,5.121745729396272) node {$v_6$};
						\draw[color=black] (0.7905042081724336,3.792307464242348) node {$e_{22}$};
						\draw [] (0.6652672701506863,5.256616278035077) circle (1.5pt);
						\draw[color=black] (0.5014958896607105,5.6226934814832585) node {$v_6$};
						\draw[color=black] (5.481375737234895E-4,4.678599641011631) node {$e_{27}$};
						\draw[color=black] (1.8116669335805224,3.704246897817612) node {$e_{21}$};
					\end{scriptsize}
				\end{tikzpicture}
				\caption{The forest graph obtained from $\vec{\Omega}$ by disconnecting every path in $U$ at its finishing vertex. } \label{dforest}
			\end{center}
		\end{figure}
		Intuitively, the depths of an edge can be viewed as the maximum length of the directed paths starting from the edge in the forest graph obtained from $\vec{\Omega}$, by disconnecting every path in $U$ at its finishing vertex. See Figure \ref{dforest} for example, where the depth of $e_{17}$ is  $\max(l_{17}+l_{16}+l_{15}, l_{17}+l_{16}+l_{14})$.  The depth of $e_{24}$ is  $\max(l_{24}+l_{18}, l_{24}+l_{17}+l_{16}+l_{15}, l_{24}+l_{17}+l_{16}+l_{14}, l_{24}+l_{21}+l_{19},l_{24}+l_{21}+l_{20})$.

		\subsection{Vertex relations}
		
	\je{In this section, we rewrite the vertex relations \rq{nc}, \rq{dc} in a way that is convenient for our proofs of shape and velocity control. The ideas here are adapted from \cite{AZ}.}

		Let $\mathbf{g}$ be a vector function, such that for any $v_i \in V \setminus \Ga,$ $g_i(t)=u_j(v_i,t),\, j  \in J^c(v_i).$  
		For $v_i \notin \Gamma$, let 
		\begin{equation} \label{hji} h_{ji}(t)=\begin{cases} {g}_i(t)+f_j(t), & j \in J^*(v_i), \\ 
				{g}_i(t), & j \in J^c(v_i). \end{cases}\end{equation}
		In what follows, we will without risk of confusion drop the subscript $i$ from $h_{ji}$.
		
		It will be useful to adapt the following notation. For edge $e_j(v_i,v_k)=(0,l_j)$, we denote $u_j^{h_j,DD,g_{k}}$ to be the solution to 
		\begin{eqnarray*}
			u_{tt}-u_{xx}+q(x)u & = & 0,\ (x,t)\in (0,l_j)\times (0,T),\\
			u(0,t)& = & h_j(t),\\
			u(l_j,t)& = & g_k(t),\\
			u(x,0)=u_t(x,0) & = & 0.
		\end{eqnarray*}
		Then it is easy to see that the solution $u$ to \rq{wave1}-\rq{ic1} satisfies $u_j(x,t)=u_j^{h_j,DD,g_{k}}(x,t),$ whose formula can be found in Section 2.3.

	Suppose first that
 $v_i \in V^* $, so $v_i$ has  only outgoing edges.  For all $j \in J^+(v_i) \cap  J^-(\Gamma)$, we have
		$u_j=u_j^{h_j,DD,0}.$
		{\bf Notation question: is the DD superscript necessary? In this paper, we have no $u^{ND}$, etc}
		For all $j \in J^+(v_i) \setminus  J^-(\Gamma)$, $u_j(x,t)=u_j^{h_j, DD, g_{k(j)}}$, where $k(j)$ is the finishing vertex of $e_j$.  Therefore, using \rq{fold2} and \rq{fold2bk}, 
		\begin{eqnarray} 
			\tilde{f}_i(t)& =&   \sum_{j \in J(v_i)}  \partial u_j(v_i,t)\nonumber \\ 
			&= &  \sum_{j \in J^+(v_i) \setminus J^-(\Gamma)} \left[ -h_j' (t) - 2 h_j'(t-2l_j) - 2 h_j'(t-4l_j) +\dots \right]\nonumber\\
			&&+ \sum_{j \in J^+(v_i) \setminus J^-(\Gamma)} \left[ 2 g_{k(j)}'(t-l_j) + 2 g_{k(j)}'(t-3l_j) +\dots \right]\nonumber\\
			&&+ \sum_{j \in J^+(v_i) \cap J^-(\Gamma)} \left[ -h_j' (t) - 2 h_j'(t-2l_j) - 2 h_j'(t-4l_j) +\dots \right]\nonumber\\
			&&+\mbox{integral terms}.\label{fw0011}
		\end{eqnarray}
		
		Now we consider $v_i\in V\setminus V^*$.
		Every $v_i \in (V \setminus V^*)\setminus \Gamma $ has incoming edges and possible outgoing edges.   For all $j \in J^-(v_i) $, $u_j(x,t)=u_j^{h_j, DD, g_i}$.  If $v_i$ has outgoing edges, for all $j \in J^+(v_i)\setminus J^-(\Gamma)$, we have $u_j(x,t)=u_j^{h_j, DD, g_{k(j)}}$;  for all $j \in J^+(v_i)\cap J^-(\Gamma)$, we have  $u_j(x,t)=u_j^{h_j, DD,0}$. Therefore
		\begin{eqnarray} \label{fw0012}
			0&
			=& \sum_{j \in J^-(v_i) \cap \je{J^+(\Gamma)}} \left[ 2  f_{j}'(t-l_j) + 2f_{j}'(t-3l_j) +\dots \right]\nonumber \\
			&&+ \sum_{j \in J^-(v_i) \cap J^+(\Gamma)} \left[ -g_i' (t) + 2 g_i'(t-2l_j) - 2 g_i'(t-4l_j) +\dots \right]\nonumber\\
			&&+\sum_{j \in J^-(v_i) \setminus \je{J^+(\Gamma)}} \left[  2 h_j'(t-l_j) + 2 h_j'(t-3l_j) +\dots \right]\nonumber\\
			&&+ \sum_{j \in J^-(v_i) \setminus J^+(\Gamma)} \left[ -g_i' (t) - 2 g_i'(t-2l_j) - 2 g_i'(t-4l_j) +\dots \right]\nonumber\\
			&&+ \sum_{j \in J^+(v_i) \setminus J^-(\Gamma)} \left[ -h_j' (t) - 2 h_j'(t-2l_j) - 2 h_j'(t-4l_j) +\dots \right]\nonumber\\
			&&+ \sum_{j \in J^+(v_i) \cap J^-(\Gamma)} \left[ -h_j' (t) - 2 h_j'(t-2l_j) - 2 h_j'(t-4l_j) +\dots \right]\nonumber\\
			&&+ \sum_{j \in J^+(v_i) \setminus J^-(\Gamma)} \left[ 2 g_{k(j)}'(t-l_j) + 2 g_{k(j)}'(t-3l_j) +\dots \right]\nonumber \\
			&&+\mbox{integral terms}.
		\end{eqnarray}
		For given control function $\mathbf{f}(t)$, each of $|V|-|\Gamma|$  equations  representing the vertex conditions \rq{nc}-\rq{dc}
		is in the form of either \eqref{fw0011} or \eqref{fw0012}.

		\

		\subsection{Proofs of shape and velocity controllability}

		\begin{thm} \label{shapevelocity}
			Let $\vec{\Omega}$ be a DAG of $\Omega$, $\{ I^*, J^*\}$ be a ST active set, and $U$ be the associated path union. Let 
			
			\begin{equation} \label{tstar}
				T_*:=\max_{j \in J} \operatorname{depth}(e_j).
			\end{equation}
			
			\begin{enumerate} 
				\item For any $\phi \in \mathcal{H}$ and $T \geq  T_*$, there exists $\textbf{f} \in \mathcal{F}_0^{T}$ on $\{I^*, J^*\}$, such that $u_j(x, T)=\phi_j(x)$, where $u_j:=u^{\textbf{f}}|_{e_j}$ and $\phi_j(x):=\phi|_{e_j}(x)$ for all $j \in J$ and $x \in [0, l_j]$.  
				\item  For any $\psi \in \mathcal{H}_0^1$ and $T > T_*$, there exists $\textbf{f} \in \mathcal{F}_2^T$ such that $(u_j)_t(x, T)=\psi_j(x)$, where $u_j:=u^{\textbf{f}}|_{e_j}$ and $\psi_j(x):=\psi|_{e_j}(x)$ for all $j \in J$ and $x \in [0, l_j]$.  
				\item  Let the greatest lower bound for both shape and velocity controllabilty on $\Omega$ be $T_{inf}$ , then $$T_{inf} =\min_{\vec{\Omega}, U}  \max_{j \in J} \operatorname{depth}(e_j, \vec{\Omega}, U). $$  
			\end{enumerate}
		\end{thm}
		{\bf Proof}
		We prove 
		{Part 2}, and thereby implicitly Part 3; the proof of Part 1 is similar. Broadly speaking, the proof reduces the solving of control functions for the graph to solving integral equations on each edge, an iterative process where the wave on an edge is related to the wave on the neighboring edges using the vertex relations developed in Section 2.3.
		
		\je{For the proof, we adopt the following notation: 
		$$H^n_*(a,b)=\{ \phi \in H^n(a,b): \ \frac{d^m\phi}{dt^m}(a)=0,\ m=0,..., n-1\} .
		$$}
		
		We show $\mathbf{f}(t) \in \mathcal{F}_2^T$ can be constructed in such a way that if the four rules listed below are followed, $\mathbf{f}(t)$ is a solution to the velocity control problem \eqref{wave1}--\eqref{ic1}.   
		
		For every $v_i \in \vec{\Omega}^-$, we pick one incoming edge $e_{m(i)}$ such that $l_{m(i)} =\min_{j \in J^-(v_i)} l_j$, and refer to $e_{m(i)}$ as the controlling edge of $v_i$. Note that for $v_i \in \Gamma \cap \vec{\Omega}^-$, its controlling edge is  the unique incoming edge of this vertex. For $v_i \in \vec{\Omega}^- \setminus \Gamma$, the choice of  controlling edge is not unique; let \begin{equation} \label{ep}
			0< \epsilon < \min_{j \in J} \{ l_j, T-T_* \}.
		\end{equation}
		The four rules for constructing the velocity control functions are:
		
		\begin{enumerate}
			\item If $v_i$ is a sink with degree $1$, let $e_j(v_k, v_i)$ be its  incoming edge, then $h_j(t)=0$ when $t$ is outside of $[T-l_j,T]$.

			\item If $v_i$ is a sink with degree more than $1$, then  $g_i(t)$ is a continuous function whose support is on $(T-\epsilon, T]$, and $g_i(T)=\phi(v_i)$.  Let $e_{m(i)}$ be its controlling edge; then $h_{m(i)}(t)=0$ when $t$ is outside of $[T-l_{m(i)}-\epsilon,T]$. For all other incoming edges $e_j$, $h_j(t)=0$ when $t$ is outside of $[T-l_j,T]$.

			\item If $v_i$ is not a sink \je{or source},  there is only one incoming edge (also denoted as $e_{m(i)}$) that is not a finishing edge of a path in $U$.  For all incoming edges $e_j(v_k, v_i)$ that are  finishing edges of paths in $U$,  $h_j(t)=0$ when $t$ is outside of $[T-l_j,T]$.  Thus from all incoming edges of $v_i$, only  $h_{m(i)}(t)$ will have an impact on $g_i(t)$ for $t\leq T$.  
			
			\item 
			
			On every $e_j(v_k, v_i)$ such that either $j\in J^*$ or $k\in I^*$, and  $v_i \in \Gamma$ 
			\begin{equation} \label{s1}
				\pa_t u_j(x,T)=u_j^{h_j', DD}(x,T)=\psi_j(x), \quad x \in [0, l_j].
			\end{equation}

			On every $e_j(v_k, v_i)$ such that either $j\in J^*$ or $k\in I^*$, and $v_i \notin \Gamma$, 
			\begin{equation} \label{s3}
				\pa_t u_j(x,T)=u_j^{h_j', DD, g_i'}(x,T)=\psi_j(x), \quad x \in [0, l_j].
			\end{equation}

		\end{enumerate}
		
		Suppose the graph contains more than one edge.  We now construct $h_j(t)$ for all $j \in J$, and derive entries of $\mathbf{f}(t)$ and $\mathbf{g}(t)$ from $h_j(t)$, following the linear ordering of $V$.  Let  $v_1, \dots, v_n$ be a linear ordering of $V$ by Lemma \ref{topoorder}.  At every $v_i$, we calculate $f_j(t)$ for $j \in J^*(v_i)$ and $\tilde{f}_i(t)$ if $i \in I^*$ based on the known $h_j(t)$ for $j  \in J^-(v_l)$ for all $v_l < v_i$.  Then we construct $h_j(t)$ for $j \in J^-(v_i)$.  The detailed steps are described below. 
		
		The first vertex, $v_1$, is a sink. 
		
		Suppose first that $v_1$ is a boundary vertex, and 
		 suppose $J^-(v_i)=\{j\}$. In what follows, 
		we identify $e_j$ with 
$(0,l_j)$ with $x=l_j$ identified with $v_1$.   We construct a function $h_j(t)$, supported on $[T-l_j, T]$, by the formula
		$$\pa_t u_j^{h_j, DD,0}(x, T)=\psi_j(x).$$
 We set $h_j(t)=0$ for $t<T-l_j$
in \rq{fold2}, and $\tilde{h}_j(t)=h_j(t+(T-l_j)).$ Then 
$$\psi_j(x)=\tilde{h}_j'(l_j-x)+\int_x^{l_j}\om_j(x,s)
\tilde{h}_j'(l_j-s)ds, \ x\in (0,l_j);
$$
here $\om_j$ is the corresponding Goursat kernel for $e_j$.
This equation is a VESK, and thus has a unique solution $\tilde{h}_j'\in H^1(0,l_j)$.
Furthermore, the equation above gives we $\tilde{h}_j'(0)=\psi_j(l_j)=0.$
Translating again in time, $h_j'\in H^1_*(T-l_j,T)$. We set $h_j'(t)=0$ for $t<T-l_j$, and integrate using $h_j(0)=0.$
We get $h_j\in H^2_*(0,T).$

		Now suppose $v_1$ is not a boundary vertex. We choose $g_1(t)$, smooth and vanishing for $t<T-\epsilon$, such  $g_1'(T)=\psi(v_1)$.  We label its controlling edge $e_{m(1)}$, which we identify with the interval $(0,l_{m(1)})$ with $v_1$ corresponding to $x=l_{m(1)}$.  By \eqref{fw0012}, \rq{fold2}, and \rq{fold2bk}, we get
		$$
			2 h_{m(1)}'(t-l_{m(1)}) 
			+2w_{m(1)}(l_{m(1)},l_{m(1)})h_{m(1)}(t-l_{m(1)})-2\int_{l_{m(1)}}^t(\om_{m(1)})_x(l_{m(1)},s)h_{m(1)}(t-s)ds
			$$
			\begin{equation} \label{fwd002}
			=   -\sum_{j \in J^-(v_1)} 
			\big ( g_1' (t)-\int_0^t(\overline{\om}_j)_x(0,s)g_1(t-s)ds
			\big ), \ t\in (T-\ep, T).  
		\end{equation}
		Here, $\om_j, \overline{\om}_j$ are the Goursat kernels associated to $e_j.$
		Denote by $\al (t)$ various functions that are known,  including the right hand side of \rq{fwd002}. Integrating this equation and using $h_{m(1)}(0)=0$, we
		can rewrite the associated integral equation as
		$$h_{m(1)}(t)+\int_0^tk(s,t)h_{m(1)}(s)\ ds=\al (t),
		\ t\in (T-l_{m(1)}-\ep, T-l_{m(1)}).
		$$
		with continuous kernel $k$.
		Thus we solve this VESK to get $h_{m(1)}(t)$ for $t \in (T-l_{m(1)}-\epsilon, T-l_{m(1)})$.   We then set 
		$h_{m(1)}(t)=0$ for $t<T-l_{m(1)}-\ep $.
		The portion of $h_{m(1)}(t)$ on $(T-l_{m(1)}, T)$ is solved as follows. Let $\tilde{h}_{m(1)}(t)=h_{m(1)}(t+(T-l_{m(1)})).$ Then \rq{s3} can be written 
\begin{eqnarray}
	\psi_{m(1)}(x)&= & 		\tilde{h}_{m(1)}'(l_{m(1)}-x)+\int_x^{T}\om_{m(1)}(x,s)\tilde{h}_{m(1)}'(l_{m(1)}-s)ds\nonumber \\
&&- \tilde{h}_{m(1)}'(-l_{m(1)}+x)-\int_{2l_{m(1)}-x}^T{\om}_{m(1)}(2l_{m(1)}-x,s)\tilde{h}_{m(1)}'(l_{m(1)}-s)ds\nonumber \\
&&			+g_1'(T-l_{m(1)}+x) +\int_{l_{m(1)}-x}^T\overline{\om}_{m(1)}(l_{m(1)}-x,s)g_1'(t-s)ds
			,\ {x\in (0,l_{m(1)}). }\nonumber
		\end{eqnarray}
		All terms in this equation are known except for the first two on the right hand side, so   we can solve for $\tilde{h}_{m(1)}'(t),$ and hence $h_{m(1)}'(t)$ for $t\in (T-l_{m(1)},T)$.
		We next show that $h_{m(1)}'(t)$ is continuous at 
		$t=T-l_{m(1)}.$
		Letting {$x\to l_{m(1)}^+$} in the equation above, we get
		$$h_{m(1)}'((T-l_{m(1)})^+)=h_{m(1)}'((T-l_{m(1)})^-),
		$$
		and hence $h_{m(1)}'$ extends to $H^1_*(T-l_{m(1)}-\epsilon,T).$
		Since the initial value of $h_{m(1)}(t)$ is $0$, one can integrate to find $h_{m(1)}(t)$ on $[ T-l_{m(1)}-\ep , T-l_{m(1)}]$.
		
		For every other incoming edge $e_j$, $j \ne m(1)$,  $h_j(t)$ is supported on $[T-l_j, T]$.  It can be computed as follows. Setting 
		$\tilde{h}_{j}(t)=h_{j}(t+(T-l_{j}))$, we have for $x\in (0,l_j)$
		\begin{equation}\label{coon2}
			\tilde{h}_{j}'(l_j-x)
			+\int_x^{l_j}\om_j(x,s)\tilde{h}_{j}'(l_j-s)ds +g_1'(T-l_j+x)
		 +\int_{l_j-x}^T\overline{\om}_{j}(l_j-x,s)g_1'(t-s)ds
			=\psi_j(x),
		\end{equation}
		we solve this VESK for $\tilde{h}_{j}'$ and hence $h_j'$.
		We then integrate $h_j'$, and extend  all $h_j(t)$ to the $t<T-l_j$ by $0$. {This completes our treatment of $v_1$. }
		
		\ 
		
		Now we make the inductive hypothesis that up to some $i \in \{2, \dots, n\}$, for every $v_l < v_i$ and all its incoming edges $e_j(v_k, v_l)$, we have constructed $h_j(t)$, such that
		
		\begin{equation} \label{abc}
			h_j(t) \in H^2_*[T-\operatorname{depth}^*(e_j), T],
		\end{equation}  
		where 
		
		\begin{equation} \label{depthsdef*}
			\operatorname{depth^*}(e_j(v_i, v_k)):=\begin{cases}
				l_j&  v_k\in \Omega^-\cap \Gamma\\
				l_j+\epsilon  & v_k \in \Omega^-\setminus \Gamma, e_j {\rm \ is \ the \ controlling \ edge \ of \ } v_i,\\
				l_j, & v_k \in \Omega^-\setminus \Gamma, e_j {\rm \ is \ not \ the \ controlling \ edge \ of \ } v_i, \\  
				l_j, & v_k \notin \Omega^-, e_j {\rm \ is  \ a \ finishing \ edge \ of \ a \ path \ in \ } U, 
				\\ l_j + \max_{r \in J^+(v_k)} \operatorname{depth}^*(e_r), & {\rm otherwise}.\end{cases} 
		\end{equation}
		
		One can prove from \eqref{depthsdef} and \eqref{depthsdef*} and by induction for any $j \in J$, 
		
		\begin{equation} \label{ddstar}
			0 \le \operatorname{depth}^*(e_j)-\operatorname{depth}(e_j) \le \epsilon.
		\end{equation}  
		Thus by \eqref{tstar}, \eqref{ep}, and \eqref{ddstar}, 
		
		\begin{equation} \label{Tg}
			T- \operatorname{depth}^*(e_j) >0 {\rm \ for \ any \ } j \in J.
		\end{equation}
		
		Moreover, as a part of the inductive hypothesis, we assume for all $v_l<v_i$ (when $i=2$, the case is trivial), we have computed $f_j(t)$ for all $j \in J^*(v_l)$ and $\tilde{f}_l(t)$ if $v_l \in V^*$, and $$f_j (t) \in \je{\hH^2} (T-\delta^*(v_l),T), \quad \tilde{f}_l(t) \in H_*^1(T-\delta^*(v_l),T), $$
		where \begin{equation} \label{dl} \delta^*(v_l):=\max_{s \in J^+(v_l)} \operatorname{depth}^*(e_s).\end{equation}  
		By \eqref{ddstar} and \eqref{dl},
		\begin{equation} \label{TgD}
			T- \delta^*(v_i) >0 {\rm \ for \ any \ } v_i\in V.
		\end{equation}
		Note that the inductive hypothesis also implies we have solved for $g_l\in H^2_*(T-depth^*(e_j),T)$ for all $v_l<v_i.$
		
		There are three cases at $v_i$: 
		
		Case 1, $v_i$ is a sink.  We compute $h_j(t)$ for all $j \in J^-(v_i)$ in the same way we compute $h_j(t)$ for all $j \in J^-(v_1)$. 
		
		Case 2,  $v_i$ is neither a sink nor a source, so $i\notin I^*.$ 
		
		We first compute $f_j(t)$ for $j \in J^*(v_i)$.  At this point we have computed $h_j(t)$ for all $j \in J^+(v_i)$, since each $j \in J^+(v_i)$ is in some $J^-(v_l)$ for some $v_l < v_i$.  For the one and only path in $U$ 
		through $v_i$, let the incoming edge of $v_i$ on $U$ be $e_{m(i)}$, and the outgoing edge of $v_i$ in $U$ be $e_{r(i)}$.   Since no control is applied on $e_{r(i)}$, \beq
		g_i(t) = h_{r(i)}(t),\ \forall t,\label{gh}
		\eeq 
		and thus $g_i$ is uniquely determined as an element of $H^2_*(T-\delta^*(v_i))$.
		
		For all $j \in J^*(v_i)$, 
		\beq \label{fhg}
		f_j(t)=h_j(t)-g_i(t).
		\eeq   
		Since each $g_s(t)=h_s(t)=0$ for $t \le T-\operatorname{depth}^*(e_s)$ for all $s \in J^+(v_i)$,  $f_j(t)=0$ for $t \le  T-\delta^*(v_i)$.   Moreover, since \je{$\psi(\cdot)$} is continuous at $v_i$, $h_j(T)=g_i(T)=\je{\psi(v_i)}$.  Thus $f_j(T)=0$. So 
		\begin{equation} \label{h01}
			f_j(t) \in \je{\hH^2} (T-\delta^*(v_i), T).
		\end{equation}
		
		Now we calculate the boundary conditions at the starting vertices of the incoming edges of $v_i$.  We first discuss $e_{m(i)}$.  Let $v_{k(m(i))}$ be its starting vertex.
		By Rule 3, \eqref{fw0012} simplifies to 
		\begin{eqnarray} 
0 &= &  2\big ( h_{m(i)}'(t-l_{m(i)}) 
			+  \om_{m(i)}(l_{m(i)},l_{m(i)})h_{m(i)}(t-l_{m(i)}) -\int_{l_{m(i)}}^t(\om_{m(i)})_x(l_{m(i)},s)h_{m(i)}(s-l_{m(i)})ds\big )\nonumber \\
&&			+\left[ 2 h_{m(i)}'(t-3l_{m(i)}) +\dots \right]\nonumber\\
&&		+ \left[ -g_i' (t) - 2 g_i'(t-2l_{m(i)}) - 2 g_i'(t-4l_{m(i)}) +\dots \right]\nonumber \\
&&		+ \sum_{j \in J^+(v_i) \setminus J^-(\Gamma)} \left[ -h_j' (t) - 2 h_j'(t-2l_j) - 2 h_j'(t-4l_j) +\dots \right]\nonumber \\
&&	+ \sum_{j \in J^+(v_i) \setminus J^-(\Gamma)} \left[ 2 g_{k(j)}'(t-l_j) + 2 g_{k(j)}'(t-3l_j) +\dots \right]\nonumber \\
&&	+ \sum_{j \in J^+(v_i) \cap J^-(\Gamma)} \left[ -h_j' (t) - 2  h_j'(t-2l_j) - 2 h_j'(t-4l_j) +\dots \right] .\label{fw004}
		\end{eqnarray}

		In \eqref{fw004}, the only unknown function is $h_{m(i)}(t)$. By the inductive hypothesis, for each
		$j\neq m(i)$, we have $h_j'(t) \in H^1_*[T-\operatorname{depth}^*(e_j), T]$, thus the sum  fourth and sixth lines of \eqref{fw004}  is in the space $H^1_*[T-\max_{j \in J^+(v_i)} \left(\operatorname{depth}^*(e_j)\right), T] $. 
		For the fifth line of \eqref{fw004},  we have $g'_{k(j)}(t) \in H^1_*[T-\operatorname{depth}^*(e_l), T]$ for some $l \in J^+(v_k(j))$. 
		Thus the entire fifth line is in the space   $H^1_*[T-\max_{j \in J^+(v_i)} \left(\operatorname{depth}^*(e_j) \right), T] $. 
		By the inductive hypotheses and \rq{gh}, $g_i'\in \je{H_*^1[T-\operatorname{depth}^*(e_{r(i)}),T]}$. Therefore the $h_{m(i)}'(t)$ obtained from \eqref{fw004} is in the following space:
		$H^1_*[T-\operatorname{depth}^*(e_{m(i)}), T-l_{m(i)}]. $  
We now use 	\eqref{fw004} to solve for $h'_{m(i)}(t)$ for $t<T-l_{m(i)}$. In what follows, we denote by $\al (t)$ various known functions. Assume first $t<3l_{m(i)}$. Recalling that
we set $h_{m(i)}(t)=0$ for $t<0$, we have
\beq 
\al (t)=  h_{m(i)}'(t-l_{m(i)}) 
			+  \om_{m(i)}(l_{m(i)},l_{m(i)})h_{m(i)}(t-l_{m(i)}) -\int_{l_{m(i)}}^t(\om_{m(i)})_x(l_{m(i)},s)h_{m(i)}(s-l_{m(i)})ds.\label{al}
\eeq
We can solve for $h_{m(i)}'(t)$ with $t<2l_{m(i)}$ similarly to how we solved \rq{fwd002}; the details are left to the reader. Next, we suppose $t<5l_{m(i)}$.
Then \rq{al} again holds, and we solve for $h_{m(i)}'$ with $t<4l_{m(i)}$. We can iterate this until we have solved for 
$h_{m(i)}'(t)$ with $t<T-l_{m(i)}$

		The portion of  $h'_{m(i)}(t)$ on the interval $(T-l_{m(i)}, T)$ controls the final velocity  on $e_{m(i)}$. Therefore, we solve 
		$$\pa_tu_{m(i)}^{h_{m(i)}, DD, g_i}(x,T)=h_{m(i)}'(T-x)+\int_x^Tw_{m(i)}(x,s)
		h'_{m(i)}(T-s)ds$$
		$$
		- h_{m(i)}'(T-2l_{m(i)}+x)+... +g_i'(T-l_{m(i)}+x)-... =\psi_{m(i)}(x).$$
		The first two terms on the right hand side are the only unknown functions here, so we
	solve this VESK for $h'_{m(i)}(t)$ on $(T-l_{m(i)}, T)$.
		{Letting $x\to l_{m(i)}^-$, we calculate 
			$h'_{m(i)}((T-l_{m(i)})^-)=h'_{m(i)}((T-l_{m(i)})^+)$, so that $h'_{m(i)}(t)$ extends to a continuous function at $x=T-l_{m(i)}$, and so $h_{m(i)}'\in H^1_*[T-\operatorname{depth}^*(e_{m(i)}),T]. $ Integrating, we get $h_{m(i)}\in H^2_*[T-\operatorname{depth}^*(e_{m(i)}),T]. $}

		For all other $j \in J^-(v_i)$, where $j \ne m(i)$, let $v_{k(j)}$ be the starting vertex of $e_j$.   We find $h_j(t)$ \je{as in \rq{coon2}} by solving 
		$$\pa_t u_j^{h_j, DD, g_i}(x,T){=h_{j}'(T-x)+\int_x^T\om_j(x,s)h_j'(T-s)ds +g_i'(T-l_{j}+x)}+... =\psi_j(x).$$  
		The terms involving $g_i$ and $\psi_j$ are all known, so we can solve this VESK for $h_j'(t)$.
		By Principle 3, $h_j(t)$ is supported on $[T-l_j, T]$,  so integrating $h_j'$ we solve for  $h_j(t)\in H_*^2[T-l_j, T]$.

		Case 3, $v_i$ is a source, so either  $i \in I^*$ \je{or
			$v_i\in \Gamma$ with its edge active.}
		
		Suppose first $i\in I^*$. For each $j \in J(v_i)$, 
		 we label $v_{k(j)}$  the vertex of $e_j$ other than $v_i$.By the inductive hypotheses, we have solved for 
		$h_j\in H_*^2(T-\operatorname{depth}^*(e_j), T)$ and $g_{k(j)}\in H^2_*(0,T).$
	Thus,  in \eqref{fw0011}, all terms on the right hand side are known.
Thus we have solved for $\tilde{f}_i\in H^1_0(T-\je{\delta^*(v_i)}, T)$. 
		
Next we solve for $f_j$ for $j\in J^*(v_i).$
		Let $e_{r(i)}$ be the one uncontrolled edge from $v_i.$
		Then we solve for $g_i(t)$ using \rq{gh}. Then $f_j\in H^2_*(T-\operatorname{depth}^*(e_j),T)$ is found using
		\rq{fhg}. By continuity of $u(*,T)$, we have $h_j(T)=g_i(T),$ and hence $f_j(T)=0.$  Similarly, by the continuity of $u_t(*,T)$, we have $h_j'(T)=g_i'(T),$ and hence $f_j(T)'=0.$ 
		Hence $f_j\in \je{H^2_0} (T-\operatorname{depth}^*(e_j),T).$
		
		\je{Finally, if $v_i\in \Gamma$ with its edge $e_j(v_i,v_k)$ active, then by the inductive hypothesis $h_i$ has already been calculated and we have $f_i=h_i\in {H}_0^2(T-\delta^*(v_i),T).$
		} This finishes Case 3.
		
		The induction has been established.$\Box$

		\begin{cor}\label{rvel}
			Assume 0 is not an eigenvalue of $A$.
			Let $\psi\in \mH^{-1}(\Omega).$
			Then there exist ${\bf f}\in \F_0^T$
			such that 
			$$u_t^{{\bf f}}(*,T)=\psi(*).
			$$
		\end{cor}
		Proof: 
Because zero is not an eigenvalue, \je{and because $\mH^{-1}$ can be constructed by the Spectral Theorem}, the mapping 
$$A:\mH^1_0(\Omega)\mapsto \mH^{-1}(\Omega)
		$$
		is an isometry.
		Thus, using \rq{ident}, {there exist $\tilde{\psi}\in \mH_0^1$ and $\tilde{\bf f}\in {\cal F}^T_2$ such that }
		\begin{eqnarray*}
			\psi & = & A\tilde{\psi}, \\
			& = & Au_t^{{\bf \tilde{f}}}(*,T), \\
			& = & \pa_t^2u_t^{{\bf \tilde{f}}}(*,T)
			= u_t^{{\bf \tilde{f}}''}(*,T),
		\end{eqnarray*}
		so it would suffice to set ${\bf f}={\bf \tilde{f}}''\in {\cal F}^T_0.$
		$\Box$
		
		\section{Exact controllability on graphs}\label{Econt}

		The proof of the exact controllability of our system uses the spectral representation of the solution $u^{\mathbf{f}}$.
		Let $\{ (\omega^2_n,\f_n): n\geq 1\}$ be the eigenvalues and normalized eigenfunctions 
		of the operator $A$, i.e. $\f_n$ solves the following eigenvalue problem:
		
		\beq \lab{sp}
		- \f_n'' + q(x) \f_n = \om_n^2 \f_n \quad \rm{on} \ \ \{\Omega\setminus V\},
		\eeq
		\beq \lab{ncf} 
		\sum_{j \in J(v_i)} \pa (\f_n)_j(v_i)=0, \  i \in I,
		\eeq
		\beq \lab{kn}
		(\f_n)_j(v_i)=(\f_n)_k(v_i), \quad  j, k \in J(v_i), i \in I.
		\eeq
		where $(\f_n)_j$ is the restriction of $\f_n$ to edge $e_j$.

		The following properties of this spectral problem are well known, see, e.g. \cite{BK}, 
		\cite{AN}.
		
		\begin{lemma}\label{specasy}
			The eigenvalues $\{ \omega_n^2\}_1^{\infty}$ associated to \rq{sp}-\rq{kn}, satisfy the relations
			$$c n \le |\omega_n| +1 \le C n,$$
			with some positive constants $c$ and $C$ independent of $n.$
			The associated unit norm eigenfunctions satisfy the relations
			$$|\f_n (v_i)| \le c_1 ,\ |(\f_n)_j'(v_i)| \le c_2 n , \quad i \in I, j \in J$$ 
			with some positive constants $c_1$, $c_2$ independent of $n,j,v_i.$
		\end{lemma} 
		
		\begin{thm} \label{full} 
			Suppose $0$ is not an eigenvalue of $A$.
			Let $\vec{\Omega}$ be a DAG of $\Omega$, and $\{I^*, J^*\}$ be a ST active set. For any $(\phi ,\psi) \in \mathcal{H}\times  (\mathcal{H}^1_0)'$, and $T >T_*$, there exists $f \in \mathcal{F}_0^{2T}$ such that $$u_j(x, 2T)=\phi_j(x),\  u_t(*, 2T)=\psi(*)$$ 
			where $u_j:=u^{f}|_{e_j}$, $\phi_j(x):=\phi|_{e_j}(x)$,  and $x \in [0, l_j]$.  
		\end{thm}

		{Proof:}
		{For this proof, we adopt the notation
		$$\ka_n^j=(\f_n)'_j (v_{\iota(j)}).
		$$}
		The solution of the IBVP \rq{wave1}--\rq{ic1} can be presented in a form of the series
		\beq \lab{ser}
		u(x,t)=\sum_{n=1}^{\infty} a_n(t)\f_n(x).
		\eeq
		Assume for the moment that the controls are all regular.
		To find the coefficients $a_n$, let $w(x,t)=\mu(t)\f_n(x)$ in the integral identity \eqref{weakwg}, where  $\mu$ is an arbitrary test function in $C^2_0[0,T].$ We obtain  the initial value problem
		$$a_n''(t)+\omega_n^2a_n(t)= \sum_{j \in J^*} \ka_n^jf_j(t) - \sum_{i \in I^*}  f'_i(t) \f_n (v_i), \quad  a_n(0)=a_n'(0)=0.
		$$
	 By variation of parameters, we get
		\beq \lab{an0}
		a_n(T)=\int_0^T \left[\sum_{j \in J^*} \ka_n^jf_j(t) - \sum_{i \in I^*}  \tilde{f}_i(t) \f_n (v_i) \right]\,\frac{\sin \om_n(T-t)}{\om_n}\,dt.
		\eeq
		We assume here that $\omega_n >0.$ If  $\omega_n =0,$ instead of $\frac{\sin \omega_n(T-t)}{\om_n}$
		we put $T-t,$ and if  $\om_n <0,$ we put $\frac{\sinh \om_n(T-t)}{\om_n}.$
		To interpret \rq{an0} for $\tilde{f}_i\in H^1(0,T)',$
		we rewrite it 
		\beq \lab{an}
		a_n(T)=\sum_{j \in J^*}\ka_n^j \int_0^T f_j(t) \frac{\sin \om_n(T-t)}{\om_n}dt- \sum_{i \in I^*} \f_n (v_i) \< \tilde{f}_i(*),\frac{\sin \om_n(T-*)}{\om_n}\> .
		\eeq
		Here and below, $\< *,*\>$ refers to either the pairing 
		$H^1(0,T)'-H^1(0,T)$ or the pairing $H^1(0,2T)'-H^1(0,2T)$;
		which one will be clear in context.
		Differentiating in $T$,
		\beq \lab{an'}
		a_n'(T)=\sum_{j \in J^*} \ka_n^j\int_0^T f_j(t) {\cos \om_n(T-t)}dt- \sum_{i \in I^*} \f_n (v_i) \< \tilde{f}_i(*),{\cos(\om_n(T-*))}\> .
		\eeq
		
		We now apply  the classical relations between the Fourier method, the moment problem, and control theory (see, e.g. \cite[Ch. III]{AI}), which  show 
		the equivalence between solvability to the moment problem \rq{an} (resp. \rq{an'}) and solvability of the shape (resp. velocity) control problem that we have already solved by  Theorem \ref{shapevelocity} and Corollary \ref{rvel}. 
		Thus
		given arbitrary sequences  $\{ a_n\},\ \{\frac{b_n}{\omega_n}\} \in \ell^2$, there exist controls $ f_{j,0},f_{j,1}\in L^2(0,T)$ and 
		$\tilde{f}_{i,0},\tilde{f}_{i,1}\in H^1(0,T)'$ such that 
		\beq \lab{an2}
		a_n=\sum_{j \in J^*}\ka_n^j \int_0^T f_{j,0}(t) \frac{\sin \om_n(T-t)}{\om_n}dt- \sum_{i \in I^*} \f_n (v_i) \< \tilde{f}_{i,0}(*),\frac{\sin \om_n(T-*)}{\om_n}\> ,
		\eeq
		\beq \lab{an2'}
		\frac{b_n}{\omega_n}=\sum_{j \in J^*}\ka_n^j \int_0^T f_{j,1}(t) \frac{\cos \om_n(T-t)}{\om_n}dt- \sum_{i \in I^*} \f_n (v_i) \< \tilde{f}_{i,1}(*),\frac{\cos \om_n(T-t)}{\om_n}\> .
		\eeq
		We now adapt an odd-even trick used in \cite{AZ} to extend
		the controls above to $[0,2T]$ and to solve the exact control problem. Specifically, we extend
		\begin{enumerate}
			\item $f_{j,0}(t)$ to be odd about $t=T$,
			\item $f_{j,1}(t)$ to be even about $t=T$.
		\end{enumerate}
		For $\tilde{f}_{i,0}$, we define its ``odd" extension $F_{i,0}$ as follows. 
		Define $R: H^1(T,2T)\mapsto H^1(0,T)$ by $\big (R\zeta \big )(T-s)=\zeta (T+s)$.
		Then
		$$\< F_{i,0},\zeta \> =\< \tilde{f}_{i,0},\zeta |_{(0,T)}\>
		-\< \tilde{f}_{i,0},R(\zeta |_{(T,2T)})\> \ ,\forall \zeta \in H^1(0,2T).$$
		For $\tilde{f}_{i,1}$, we define its ``even" extension $F_{i,1}$ by
		$$\< F_{i,1},\zeta \> =\< \tilde{f}_{i,1},\zeta|_{(0,T)}\>
		+\< \tilde{f}_{i,1},R(\zeta |_{(T,2T)})\> ,\forall \zeta \in H^1(0,2T) .$$
		Then we define
		$$
		f_j(t)=\frac{f_{j,0}(t)+f_{j,1}(t)}{2}\in L^2(0,2T),
		{F}_i=\frac{F_{i,0}+F_{i,1}}{2} \in H^1(0,2T)'.
		$$
		Applying these to \rq{an2},\rq{an2'}, we get
		\beq \lab{an3}
		a_n=\sum_{j \in J^*} \ka_n^j\int_0^{2T} f_{j}(t) \frac{\sin \om_n(T-t)}{\om_n}dt- \sum_{i \in I^*} \f_n (v_i) \< F_{i}(*),\frac{\sin \om_n(T-*)}{\om_n}\> ,
		\eeq
		\beq \lab{an3'}
		\frac{b_n}{\omega_n}=\sum_{j \in J^*} \ka_n^j\int_0^{2T} f_{j}(t) \frac{\cos \om_n(T-*)}{\om_n}dt- \sum_{i \in I^*} \f_n (v_i) \< F_{i}(*),\frac{\cos \om_n(T-*)}{\om_n}\> .
		\eeq
These equations are equivalent to 
		\beq \lab{an4}
		\frac{b_n}{\omega_n}+ia_n=\sum_{j \in J^*}\ka_n^j \int_0^{2T} f_{j}(t) \frac{e^{i \om_n(T-t)}}{\om_n}dt- \sum_{i \in I^*} \f_n (v_i) \< F_{i}(*),\frac{e^{i \om_n(T-t)}}{\om_n}\> .
		\eeq
		\beq \lab{an4'}
		\frac{b_n}{\omega_n}-ia_n=\sum_{j \in J^*} \ka_n^j\int_0^{2T} f_{j}(t) \frac{e^{-i \om_n(T-t)}}{\om_n}dt- \sum_{i \in I^*} \f_n (v_i) \< F_{i}(*),\frac{e^{-i \om_n(T-t)}}{\om_n}\> .
		\eeq
		Solvability of these moment problems is clearly equivalent to solvability of 
		\beq \lab{an5}
		\frac{b_n}{\omega_n}+ia_n=\sum_{j \in J^*} \ka_n^j\int_0^{2T} f_{j}(t) \frac{e^{i \om_n(2T-t)}}{\om_n}dt- \sum_{i \in I^*} \f_n (v_i) \< F_{i}(*),\frac{e^{i \om_n(2T-t)}}{\om_n}\> ,
		\eeq
		\beq \lab{an5'}
		\frac{b_n}{\omega_n}-ia_n=\sum_{j \in J^*} \ka_n^j\int_0^{2T} f_{j}(t) \frac{e^{-i \om_n(2T-t)}}{\om_n}dt- \sum_{i \in I^*} \f_n (v_i) \< F_{i}(*),\frac{e^{-i \om_n(2T-t)}}{\om_n}\> .
		\eeq
	and hence equivalent to 
		\beq \lab{an6}
		\frac{b_n}{\omega_n}=\sum_{j \in J^*} \ka_n^j\int_0^{2T} f_{j}(t) \frac{\cos({ \om_n(2T-t)})}{\om_n}dt- \sum_{i \in I^*} \f_n (v_i) \< F_{i}(*),\frac{\cos({ \om_n(2T-t)})}{\om_n}\> .
		\eeq
		\beq \lab{an6'}
		a_n=\sum_{j \in J^*}\ka_n^j \int_0^{2T} f_{j}(t) \frac{\sin({ \om_n(2T-t)}}{\om_n})dt- \sum_{i \in I^*} \f_n (v_i) \< F_{i}(*),\frac{\sin({ \om_n(2T-t)}}{\om_n})\> .
		\eeq

		Thus, for arbitrary sequences $\{ a_n\},\{ b_n/\om_n\} \in \ell^2$, there exist functions 
		$f_j\in L^2(0,2T), F_i\in H^1(0,2T)'$ solving the moment problems \rq{an6}\rq{an6'}. Inspection of \rq{an2},\rq{an2'}, and using  the classical relations between the Fourier method, the moment problem, and control theory (see, e.g. \cite[Ch. III]{AI}), shows that this is equivalent to solving the exact control problem. Theorem \ref{full} is proved.

		\noindent {\bf  Acknowledgments}\\
		The research of Sergei Avdonin was  supported  in part by the National Science Foundation,
		grant DMS 1909869.

	\end{document}